\documentclass{elsart}

\usepackage{amsmath,amssymb,theorem}
\usepackage{psfrag}
\usepackage{graphicx}

{\theorembodyfont{\rmfamily} 
  }
{\theorembodyfont{\rmfamily} 
  \newtheorem{example}{Example}}  
\newtheorem{mynote}{Note}

\begin{document} 

%\mycaptionstyle
%\captionstyle{fancy}

%%\RCS $RCSfile: paper.tex,v $
%%\RCS $Revision: 1.6 $
%%\RCS $Date: 2006/01/04 20:43:14 $

\newcommand{\en}[1]{(\ref{eq:#1})}
\newcommand{\leavethisout}[1]{}
\newcommand{\psfragfontsize}{\footnotesize}
\newcommand{\ds}[1]{{\displaystyle #1}}
\newcommand{\textsub}[2]{{#1}_{\text{{\em #2}}}}
\newcommand{\order}[1]{{\cal O}\left(#1\right)}
\newcommand{\tstar}{t^\ast}
\newcommand{\xstar}{x^\ast}
\newcommand{\figfont}[1]{{\small #1}}
\newcommand{\sech}{\operatorname{sech}}

\renewcommand{\order}[1]{{O}\!\left(#1\right)}

\runauthor{Soheili and Stockie}
\begin{frontmatter}
  \title{A moving mesh method \\
    with variable mesh relaxation time}
  %%\thanks[revision]{
  %% Latest Revision: \RCSRevision\ (\RCSDate).  Printed: \today.\\
  %%
  %%
  \author[ars]{Ali Reza Soheili}
  and
  \author[jms]{John M. Stockie\thanksref{fund}}
  \thanks[fund]{This author was supported by a grant from the
    Natural Sciences and Engineering Research Council of Canada
    (NSERC).\\
    \emph{E-mail addresses:} soheili@math.usb.ac.ir (A. R. Soheili),
    stockie@math.sfu.ca (J. M. Stockie).}
  \address[ars]{Department of Mathematics, University of Sistan and
    Baluchestan, Zahedan, Iran}
  \address[jms]{Department of Mathematics, Simon Fraser University,
    Burnaby, BC, Canada}
  
  %%\date{\today}
  %%\maketitle
  
  \begin{abstract}
    We propose a moving mesh adaptive approach for
    solving time-dependent partial differential equations.  The motion
    of spatial grid points is governed by a moving mesh PDE (MMPDE) in
    which a \emph{mesh relaxation time} $\tau$ is employed as a
    regularization parameter.  Previously reported results on MMPDEs
    have invariably employed a constant value of the parameter $\tau$.
    We extend this standard approach by incorporating a variable
    relaxation time that is calculated adaptively alongside the solution
    in order to regularize the mesh appropriately throughout a
    computation.  We focus on singular problems involving self-similar
    blow-up to demonstrate the advantages of using a variable relaxation
    time over a fixed one in terms of accuracy, stability and
    efficiency.  

    \noindent\emph{AMS Classification:}  
    65M50,
    65M06,
    35K57
  \end{abstract}

  \begin{keyword}
    Moving mesh method;
    Self-similar blow-up;
    Relaxation time
    %% Temporal smoothing
  \end{keyword}

\end{frontmatter}

%%%%%%%%%%%%%%%%%%%%%%%%%%%%%%%%%%%%%%%%%%%%%%%%%%%%%%%%%%%%%%%%%%%%%%%%
\section{Introduction}
\label{sec:intro}

Moving mesh methods have been employed widely to approximate solutions
of partial differential equations which exhibit large solution
variations, such as shock waves and boundary or interior layers.
Several moving mesh approaches have been derived and many authors have
discussed the significant improvements in accuracy and efficiency that
can be achieved with respect to fixed mesh methods
\cite{cm98,hrr94a,hrr94b,ls06,sc86,smr01,ttz03}.

The moving mesh PDE (or MMPDE) approach has proven particularly
effective in solving nonlinear PDEs that exhibit solutions having some
type of singularity, such as self-similar blow-up~\cite{bhr96} or moving
fronts~\cite{bhjj06,smr01}.  For blow-up problems in particular, moving
mesh methods permit a detailed study of the singularity formation with a
degree of accuracy and efficiency that is simply not possible using
fixed mesh methods.  The primary advantage of the moving mesh approach
stems from its ability to exploit special features of the solution (such
as self-similarity) and build them directly into the numerical scheme.

In the MMPDE approach, a separate PDE is derived to evolve the mesh
points in such a way that they tend towards an equidistributed mesh at
steady state, in the sense that the mesh points are positioned in space
so as to equally distribute some measure of the solution error.  The
MMPDE is coupled nonlinearly to the physical PDE of interest, and both
PDEs are solved simultaneously.  A key parameter in the moving mesh
equation is the \emph{mesh relaxation time,} usually denoted as $\tau$;
the exact equidistribution equation is notoriously
ill-conditioned~\cite{asc99,pet87,lpr99} and so $\tau$ acts to
regularize the mesh evolution in time.  The philosophy behind
introducing temporal smoothing, instead of equidistributing exactly, is
that the mesh need not be solved to the same level of accuracy as the
physical PDE; in fact, solution accuracy can still be significantly
improved over fixed mesh methods by only approximately equidistributing
the mesh.

In previous results reported in the literature, the mesh relaxation time
is invariably taken to be a constant for any given simulation.
Furthermore, Huang, Ren and Russell observed in \cite{hrr94b} that
\emph{``while the parameter $\tau$ is critical, in our experience the
  numerical methods are relatively insensitive to the actual choice of
  $\tau$ in applications,''} and similar comments were made in
\cite{bhr96,hrr94a}.  However, it is essential to keep in mind that
these observations were made for problems in which the range of time
scales present in the solution was fairly limited.  In practice, $\tau$
must be tuned manually to optimize the behaviour of the computed mesh,
and sometimes even to obtain a convergent numerical solution.

The main purpose of this paper is to consider situations where taking
constant $\tau$ may not be appropriate.  Keeping in mind that $\tau$ can
be interpreted as a time scale for the mesh motion, then $\tau$ should
in fact be taken as a solution-dependent parameter, because as
singularities form, intensify, propagate, and dissipate, the speed of
solution variations (and hence also of the mesh points) in a given
computation may vary a great deal.  By no means are we suggesting that a
variable $\tau$ is necessary in all moving mesh calculations.
Nonetheless, there is some advantage to be gained by having an algorithm
that is capable of determining the value of $\tau$ automatically as part
of the solution process without requiring the user to determine its value
through trial and error (since the complicated nonlinear coupling
between solution and mesh in the MMPDE approach means there is no way to
know the value of $\tau$ \emph{a priori}).  The main purpose of this
paper is to demonstrate, by means of specific examples, that varying the
mesh relaxation parameter throughout a computation can be of significant
advantage in terms of both accuracy and efficiency.  We will present an
approach for adaptively selecting $\tau$ in such a way that the temporal
evolution of the mesh is optimal in an appropriate sense.

This paper is organized as follows.  In Section \ref{sec:mmpde}, we
briefly review moving mesh methods in which the mesh equation
incorporates a relaxation time $\tau$.  The main motivating example for
introducing an adaptive strategy for choosing a time-dependent mesh
smoothing parameter comes from a class of nonlinear parabolic equations
exhibiting self-similar blow-up behaviour; we therefore introduce in
Section \ref{sec:blowup} the blow-up model equation, and motivate a
particular choice of $\tau(t)$ which is suggested by the analysis of
blow-up problems.  Numerical experiments are then presented in Section
\ref{sec:results} to illustrate the advantages of this modified moving
mesh approach in terms of both accuracy and efficiency.

%%%%%%%%%%%%%%%%%%%%%%%%%%%%%%%%%%%%%%%%%%%%%%%%%%%%%%%%%%%%%%%%%%%%%%%%
\section{The moving mesh method}
\label{sec:mmpde}

The evolution of a moving computational grid can be viewed as a
discretization of a one-to-one, time-dependent coordinate mapping.  Let
$x$ and $\xi$ denote the physical and computational coordinates
respectively, and define a coordinate transformation by
\begin{gather*}
  x = x(\xi,t) \quad \text{where} \quad x(0,t) = 0 \quad \text{and} \quad
  x(1, t)=1,
\end{gather*}
where both $x$ and $\xi$ are assumed to lie in interval $[0, 1]$.  The
computational coordinate is discretized on a uniform mesh given by
$\xi_i = i/N$, where $i=0,1,2,\dots,N$ and $N$ is a positive
integer.  The corresponding non-uniform mesh is denoted by
\begin{align*}
  0 = x_0 < x_1(t) < x_2(t)< \cdots < x_{N-1}(t) < x_N = 1.  
\end{align*}

A key ingredient of the moving mesh approach is the \emph{monitor
  function}, $M(x,t$), which is chosen to be some approximate measure of
the solution error.  For a given monitor function, the mesh point locations
$x_i(t)$ could be required to satisfy the following equidistribution
principle (EP) for all values of time $t$ \cite{hrr94b}:
\begin{gather*} 
  \int_{x_{i-1}(t)}^{x_i(t)} M (x, t) d x = \frac{1}{N}
  \int_0^1 M(x, t) d x = \frac{\theta(t)}{N}, 
\end{gather*} 
or equivalently
\begin{gather} 
  \int_{0}^{x_i(t)} M (x^\prime, t) \,dx^\prime =  \frac{i}{N}
  \theta(t)=\xi_i \, \theta(t),  
  \label{eq:ep}
\end{gather} 
where $\theta(t) = \int_0^1 M(x^\prime, t)\, dx^\prime$.  This EP is
intended to concentrate points in regions where $M$ (and hence also the
solution error measure) is large, thereby placing fewer points in areas
where the error is small.  Differentiating \en{ep} yields an
equivalent differential form
\begin{gather} 
  \frac{\partial}{\partial \xi} \left( M
    \frac{\partial x}{\partial \xi}\right)(\xi, t) = 0,
  \label{eq:mmesh}
\end{gather} 
where $x(0,t)=0$ and $x(1,t)=1$.

Solving the elliptic equation \en{mmesh} directly is often problematic,
since it introduces a nonlinear coupling between the mesh and the
solution through the dependence of $M$ on the solution $u$.
Furthermore, when the physical and mesh PDEs are discretized, they take
the form of an index-2 DAE system which is very stiff in practice and
also typically ill-conditioned~\cite{asc99,pet87,lpr99}.  As a result,
it is usually much more attractive to relax the requirement of exact
equidistribution by introducing a relaxation time $\tau$ into the
problem.  A dynamic moving mesh equation can be derived by requiring the
mesh to satisfy the above EP at a later time $t+\tau$ instead of at $t$.
Because $\theta(t)$ has been eliminated from the differential form of
the EP, the mesh must satisfy
\begin{gather}
  \frac{\partial}{\partial
    \xi}\left[ M\left( x(\xi, t+\tau), t+\tau\right)\,\frac{\partial
    }{\partial \xi}x(\xi, t+\tau)\right] = 0.
  \label{eq:mmesh-tau}
\end{gather}
By expanding the terms $\frac{\partial }{\partial \xi}x(\xi, t+\tau)$
and \mbox{$M(x(\xi, t+\tau), t+\tau)$} in Taylor series and dropping
certain higher order terms, a variety of different MMPDEs can be derived
\cite{hrr94b}.  In this paper we will employ two specific moving mesh
equations: 
\begin{align} 
  \text{MMPDE4:} & &
  \tau \frac{\partial }{\partial \xi} 
  \left(M\frac{\partial \dot{x} }{\partial \xi}\right) &= 
  - \frac{\partial }{\partial
    \xi}\left(M \frac{\partial x
    }{\partial \xi}\right), \label{eq:mmpde4} \\
  \text{MMPDE6:} & &
  \tau \frac{\partial ^2 \dot{x} }{\partial \xi^2} &= - \frac{\partial
  }{\partial \xi}(M \frac{\partial x}{\partial \xi}).  \hspace{2cm}
  \label{eq:mmpde6} 
\end{align}
The relaxation parameter $\tau$ can also be thought of as introducing
\emph{temporal smoothing} into the mesh.  The equivalent derivation for
a time-dependent $\tau(t)$ is given in Appendix~\ref{sec:taylor}.

The choice of monitor function $M$ in this method is somewhat arbitrary;
in general, $M$ can be any given function of $u$, or it may also be
based on some error estimate determined numerically based on discrete
solution values.  A commonly used monitor function is
$M=\sqrt{1+u_x^2}$, which equidistributes the arclength of the solution
$u$.  However, this choice of $M$ often behaves very badly in
simulations, concentrating too many points in singularities and making
the coupled problem excessively stiff~\cite{lpr99,smr01}.  Other common
monitor function are chosen either for analytical reasons (such as the
form $M=|u|^p$ for self-similar blow-up problems \cite{bhr96}) or for
practical considerations (such as the component-averaged monitor
developed in \cite{smr01} for hyperbolic systems).

In this work, we discretize the mesh equation using centered finite
differences in space, which yields for MMPDE6: 
\begin{subequations}\label{eq:discrete-mmpde}
\begin{gather}
  %%M_{i+1/2}\left( \dot{x}_{i+1} - \dot{x}_i
  %%\right) - M_{i-1/2} \left( \dot{x}_{i} -
  %%\dot{x}_{i-1} \right) = - \frac{E_i}{\tau} 
  %%\label{eq:discrete-mmpde4} , 
  %%\\
  \dot{x}_{i+1} - 2 \dot{x}_i + \dot{x}_{i-1} = - \frac{E_i}{\tau} .
  \label{eq:discrete-mmpde6}  
\end{gather}
The quantity $E_i$ represents a centered approximation to the term on
the right hand side of MMPDE6 given by
\begin{gather}
  E_i = M_{i+1/2}\left( x_{i+1} - x_i \right) -
  M_{i-1/2}\left( x_{i} - x_{i-1} \right) ,
  \label{eq:discrete-rhs}
\end{gather}
\end{subequations}
where $M_{i+1/2} = \half(M_i + M_{i+1})$, $M_i=M(u_i)$, and $u_i\approx
u(x_i,t)$ is an approximation of the solution at grid point $x_i$.  The
discretization for MMPDE4, which is also employed here, is carried out
in a similar fashion.  It turns out to be very important to smooth the
monitor function in space as well as in time in order to avoid
oscillatory errors in mesh locations which can then feed into the
solution.  To this end, the discrete monitor function values $M_i$ are
usually replaced with smoothed versions
\begin{gather} 
  \widetilde{M_i} =\left( \frac{\sum_{j=i-ip}^{i+ip} M_j^2
      \left(\frac{\gamma}{1+\gamma}\right)^{|j-i|}}{\sum_{j=i-ip}^{i+ip}
      \left(\frac{\gamma}{1+\gamma}\right)^{|j-i|}}\right)^{1/2},
      \label{eq:mtilde} 
\end{gather}
where $\gamma$ and $ip$ are parameters that must be chosen
appropriately.

\subsection{The effect of $\tau$ on mesh movement}
\label{sec:tau-ex}

Because of the importance of temporal smoothing in the present work, it
is helpful to first consider some illustrative examples that elucidate
the effect of $\tau$ on the computed mesh motion.  For this purpose, we
consider first three examples wherein $u(x,t)$ is a given function, so
that the mesh equation is uncoupled from the physical PDE.

%% EXAMPLE 1.
\begin{example}
  We first consider 
  \begin{gather} 
    u(x,t) = e^{-10 \pi^2t} \sin(\pi x),
    \label{eq:ugiven1}
  \end{gather}
  for values of $x\in [0,1]$ and $t\in[0,10]$, which was used in
  \cite{cfl86} to study the stability of various moving mesh equations,
  and also as a numerical example in \cite{hrr94b}.  The mesh equation
  is chosen to be MMPDE6, which is discretized using standard centered
  finite differences. For the purposes of this example, we use the
  arclength monitor function $M=\sqrt{1+u_x^2}$ with spatial smoothing
  parameters $\gamma=2$ and $ip=4$.  The spatial domain is divided into
  $N=100$ mesh points and the value of $\tau$ is taken to be a constant
  ranging from $10^0$ down to $10^{-5}$.  The mesh points are initially
  uniformly distributed, so that the mesh undergoes a rapid initial
  transient as the grid points are driven towards equidistribution by
  the MMPDE.  The speed of this initial transient is governed by the
  choice of the mesh relaxation time parameter.  Also, since
  $u_x(x,t)\rightarrow 0$ in the limit as $t\rightarrow\infty$, then for
  the arclength monitor $M\rightarrow 1$ as $t\rightarrow +\infty$;
  therefore, the equidistributed mesh should tend over long times to a
  uniform mesh in space.
  
  Figure \ref{fig1u} shows solution curves and Figure \ref{fig1} the
  mesh trajectories (i.e., contours of 
  $\xi(x,t)$) for the above example using MMPDE6 and a uniform initial
  mesh.  These results demonstrate a few important points.  First, the
  ability of the moving mesh to respond to rapid solution transients (in
  this case represented by the initial transient mesh redistribution) is
  governed in large part by the choice of $\tau$.  In particular, if
  $\tau$ is taken too large, then the mesh is incapable of adapting
  sufficiently well to keep up with the solution, which is 
  easily seen here in the case $\tau = 10^{-1}$.  Secondly, once
  $\tau$ is taken small enough, there is no longer any significant
  change in the mesh locations, as seen by comparing the mesh
  trajectories when $\tau=10^{-3}$ and $\tau=10^{-5}$ (ignoring the
  initial transients which are not physical but driven solely by the
  artificially chosen initial uniform mesh).  When $\tau$ is taken
  smaller than $10^{-5}$ there is no visible change in either the
  computed mesh or the time stepping behaviour.
  \begin{figure}[htbp]
    \centering
    \includegraphics[width=0.45\textwidth]{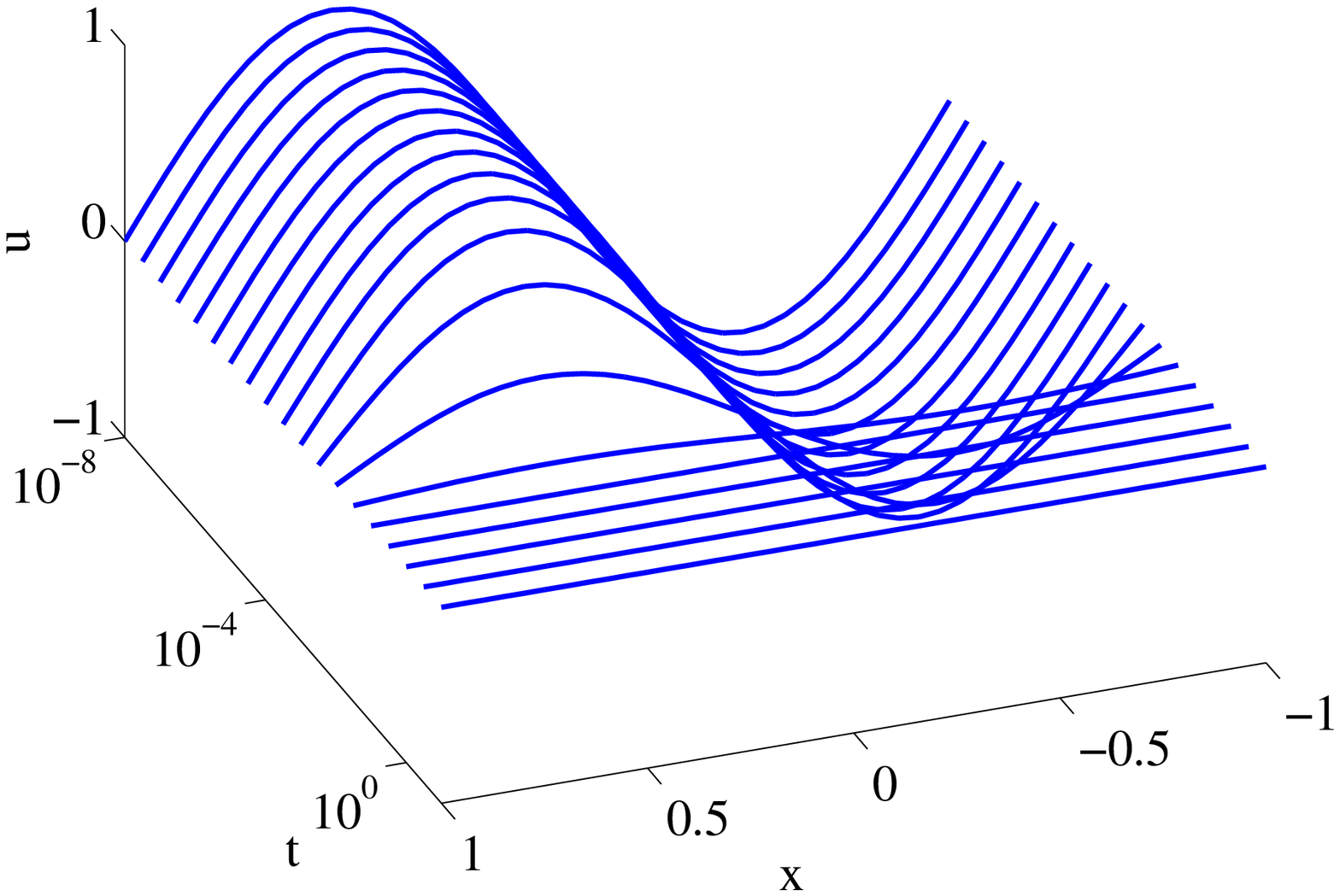}
    \quad 
    \includegraphics[width=0.45\textwidth]{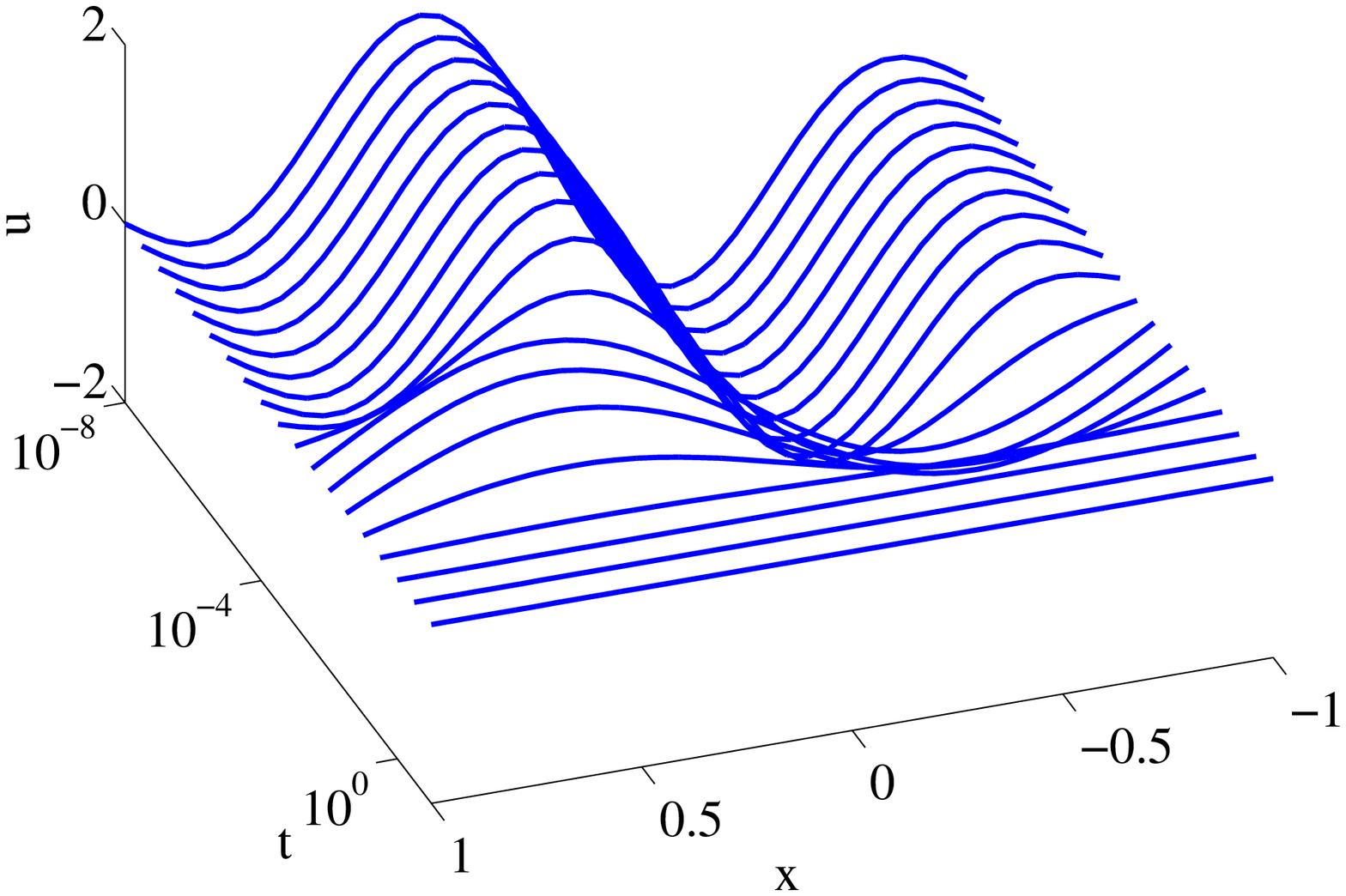}
    \caption{Solution profiles for the function \en{ugiven1} used in
      Example 1 (left) and \en{ugiven2} from Example 2 (right).} 
    \label{fig1u}
  \end{figure}
\begin{figure}[hbtp]
    \centering
    \psfragfontsize
    \psfrag{x}[c][b]{$x$}
    \psfrag{t}[Bc][c]{$t$}
    \psfrag{ t}[Bc][c]{$\Delta t$}
    \psfrag{D}{}
    \psfrag{Mesh trajectories}{}
    \psfrag{Time step}{}
    \begin{tabular}{ccc}
      \figfont{$\tau = 10^{-1}$, $\text{CPU}=7\;s$}  &
      \figfont{$\tau = 10^{-3}$, $\text{CPU}=12\;s$} &
      \figfont{$\tau = 10^{-5}$, $\text{CPU}=13\;s$} \\
      \includegraphics[width=0.3\textwidth]{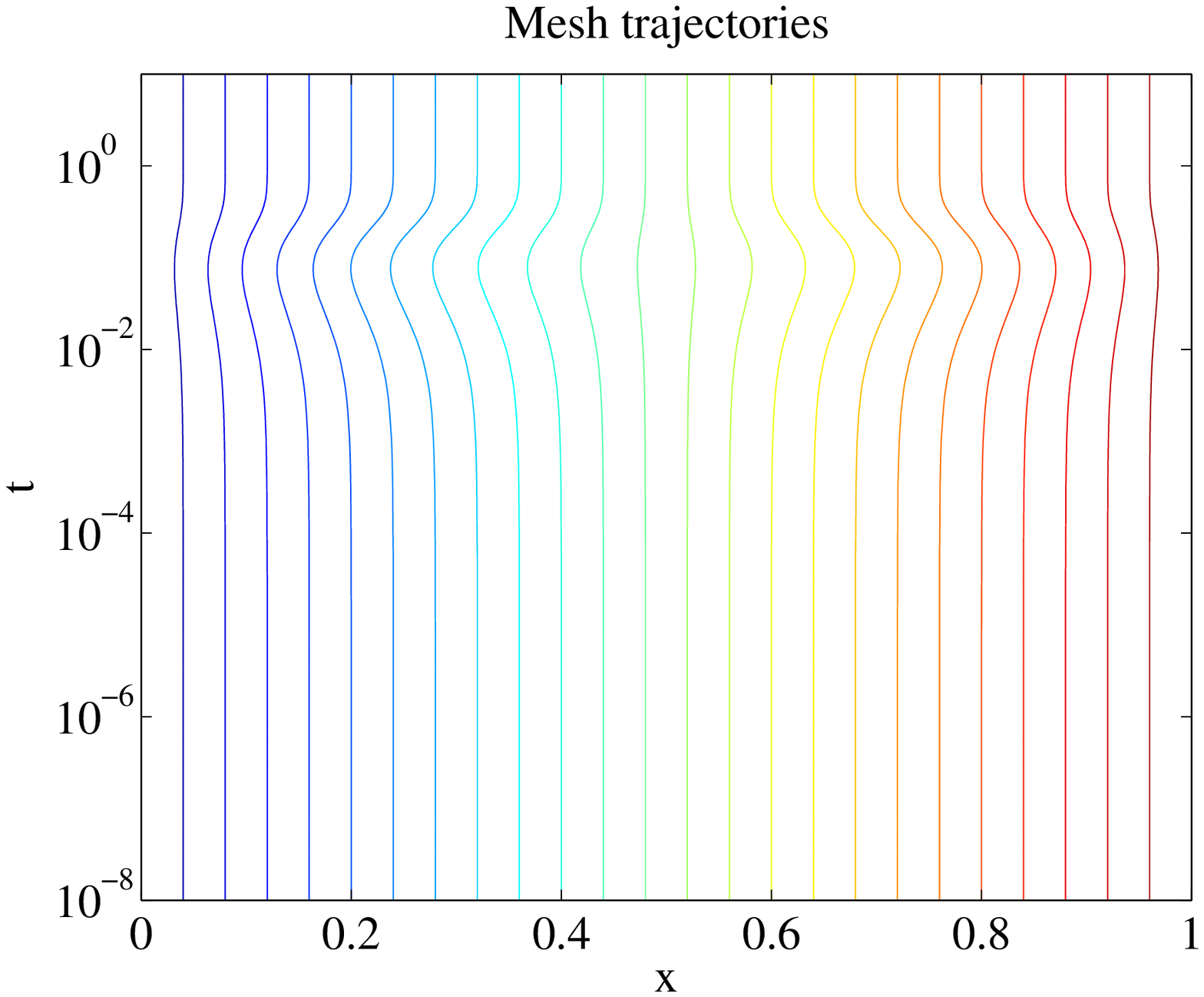} &
      \includegraphics[width=0.3\textwidth]{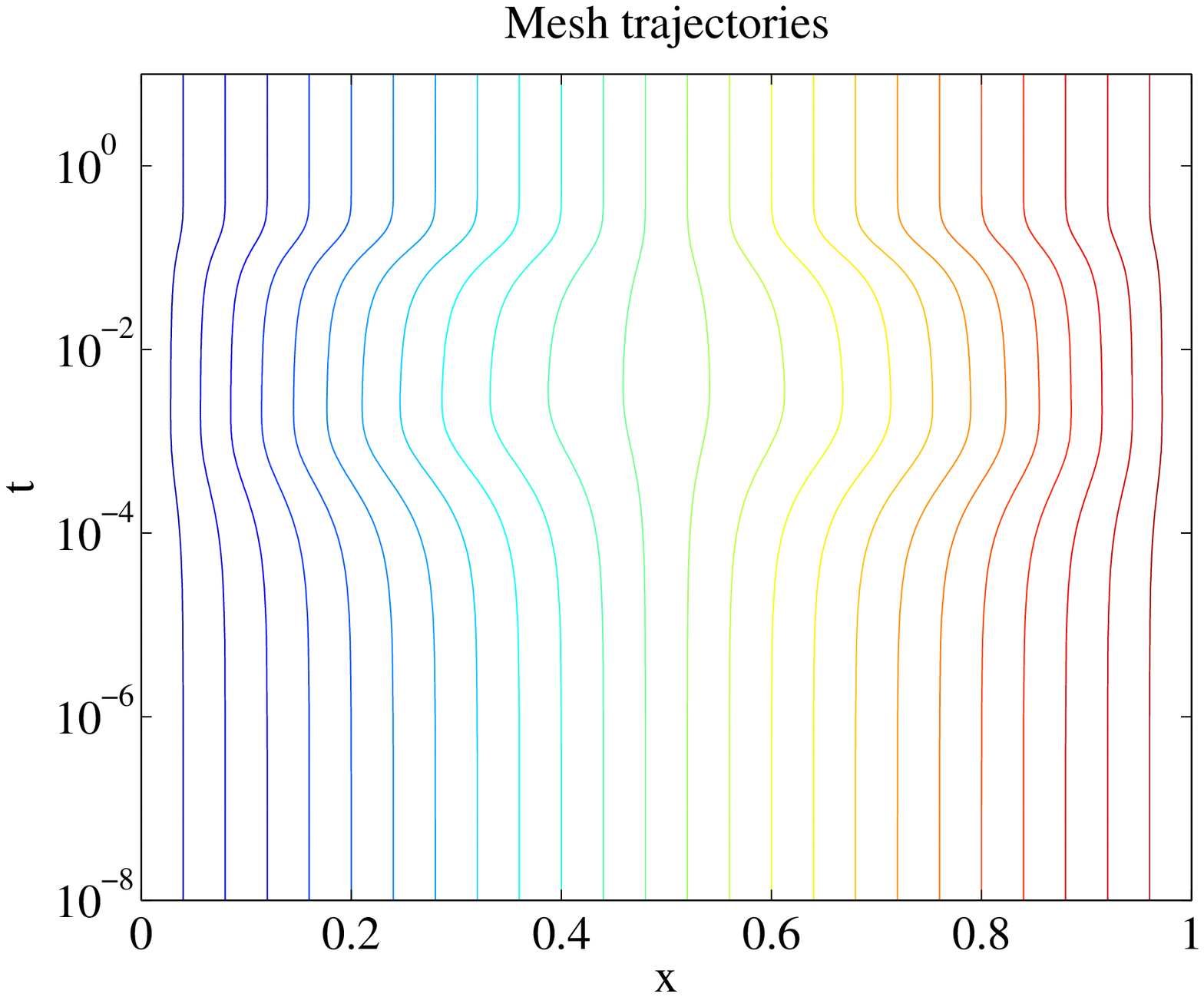} &
      \includegraphics[width=0.3\textwidth]{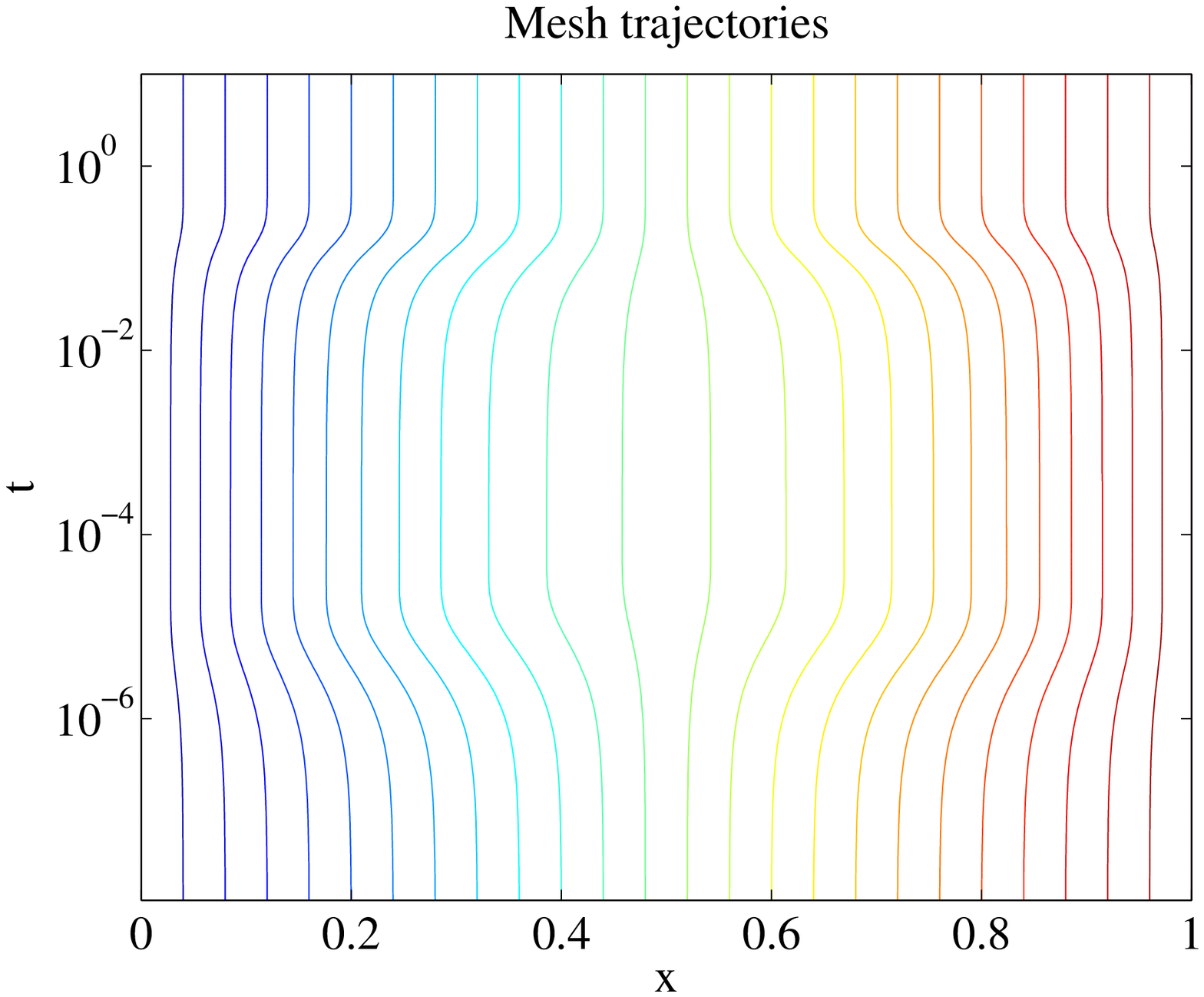}\\
      \psfrag{t}[c][b]{$t$}
      \includegraphics[width=0.3\textwidth]{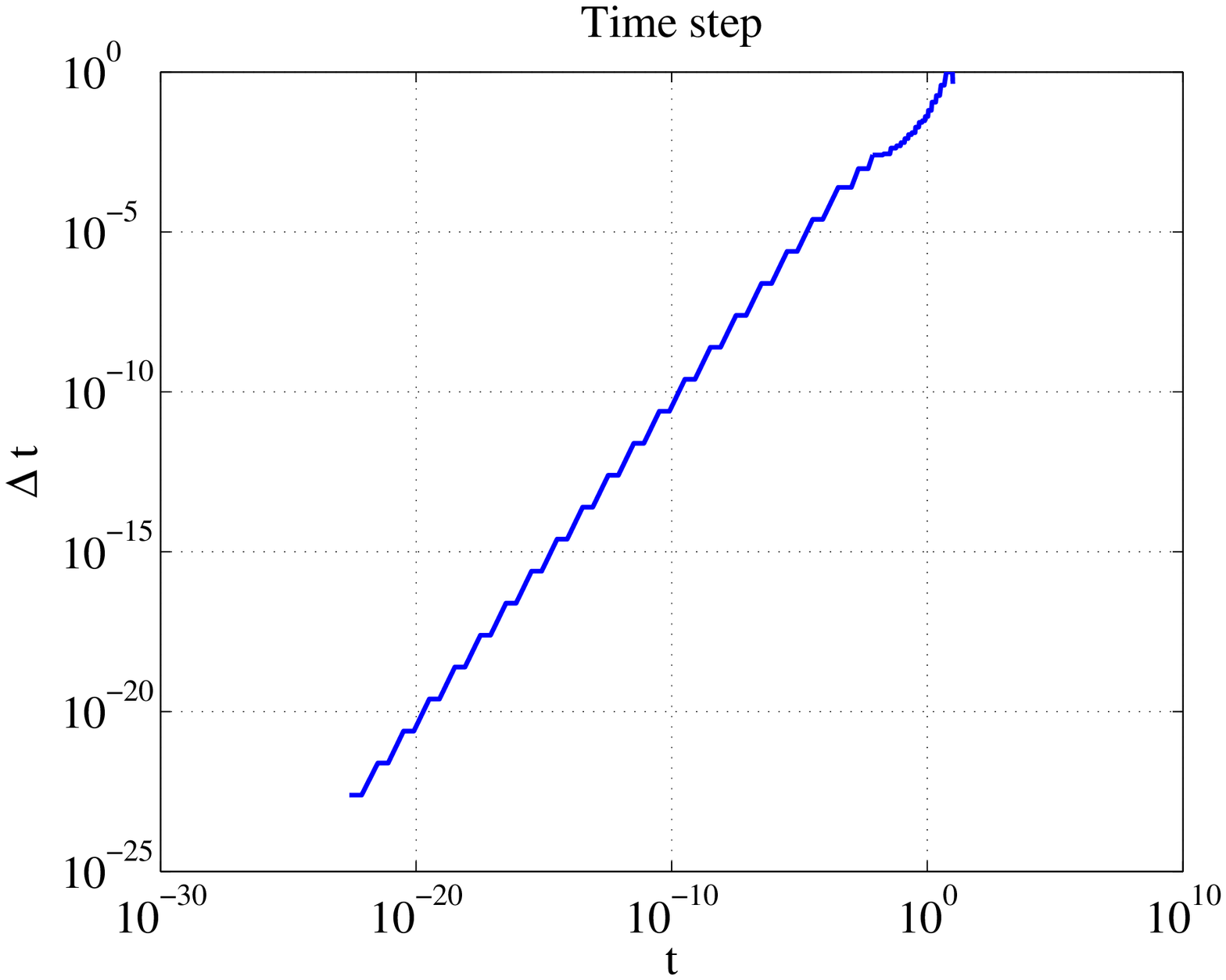} &
      \psfrag{t}[c][b]{$t$}
      \includegraphics[width=0.3\textwidth]{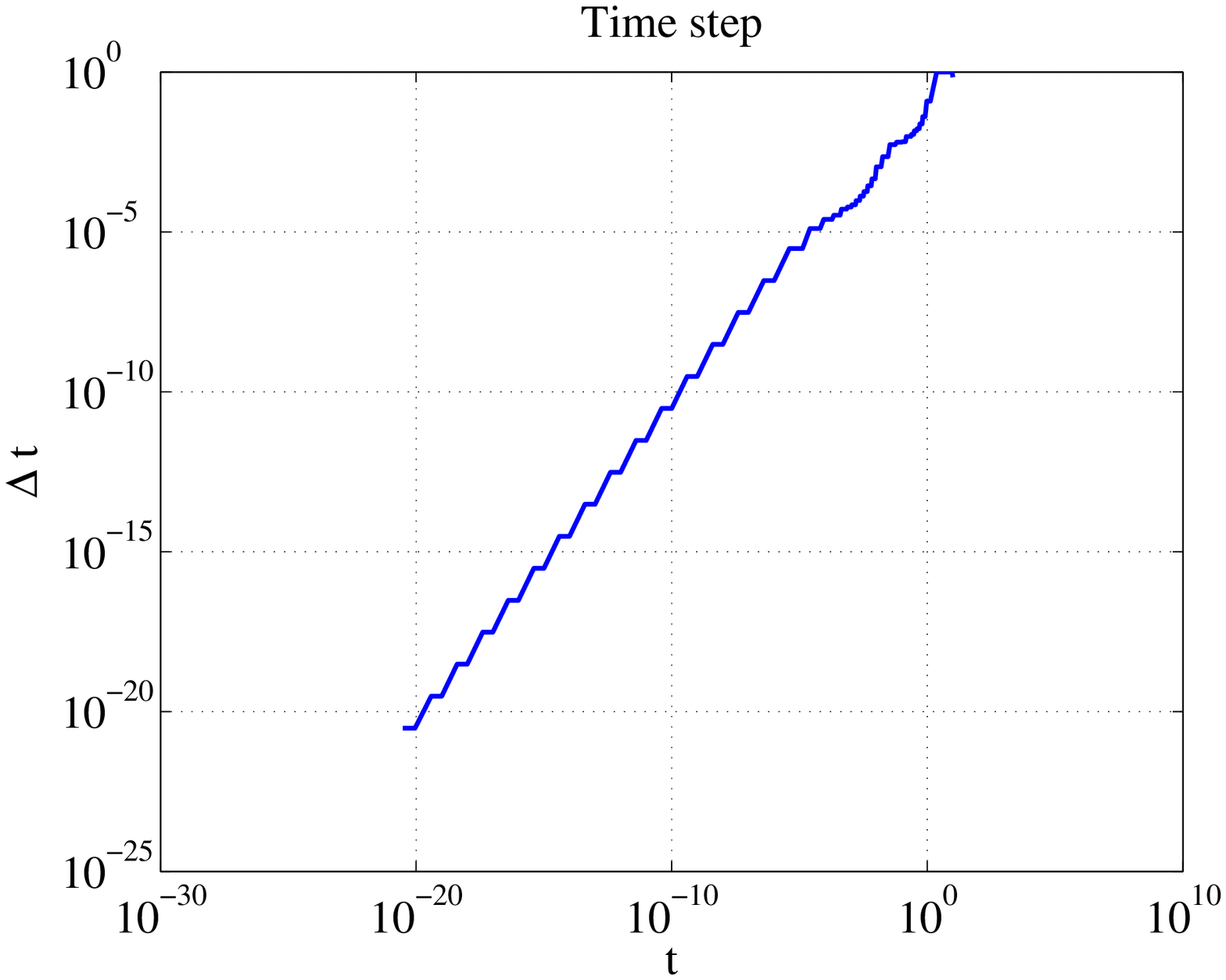} &
      \psfrag{t}[c][b]{$t$}
      \includegraphics[width=0.3\textwidth]{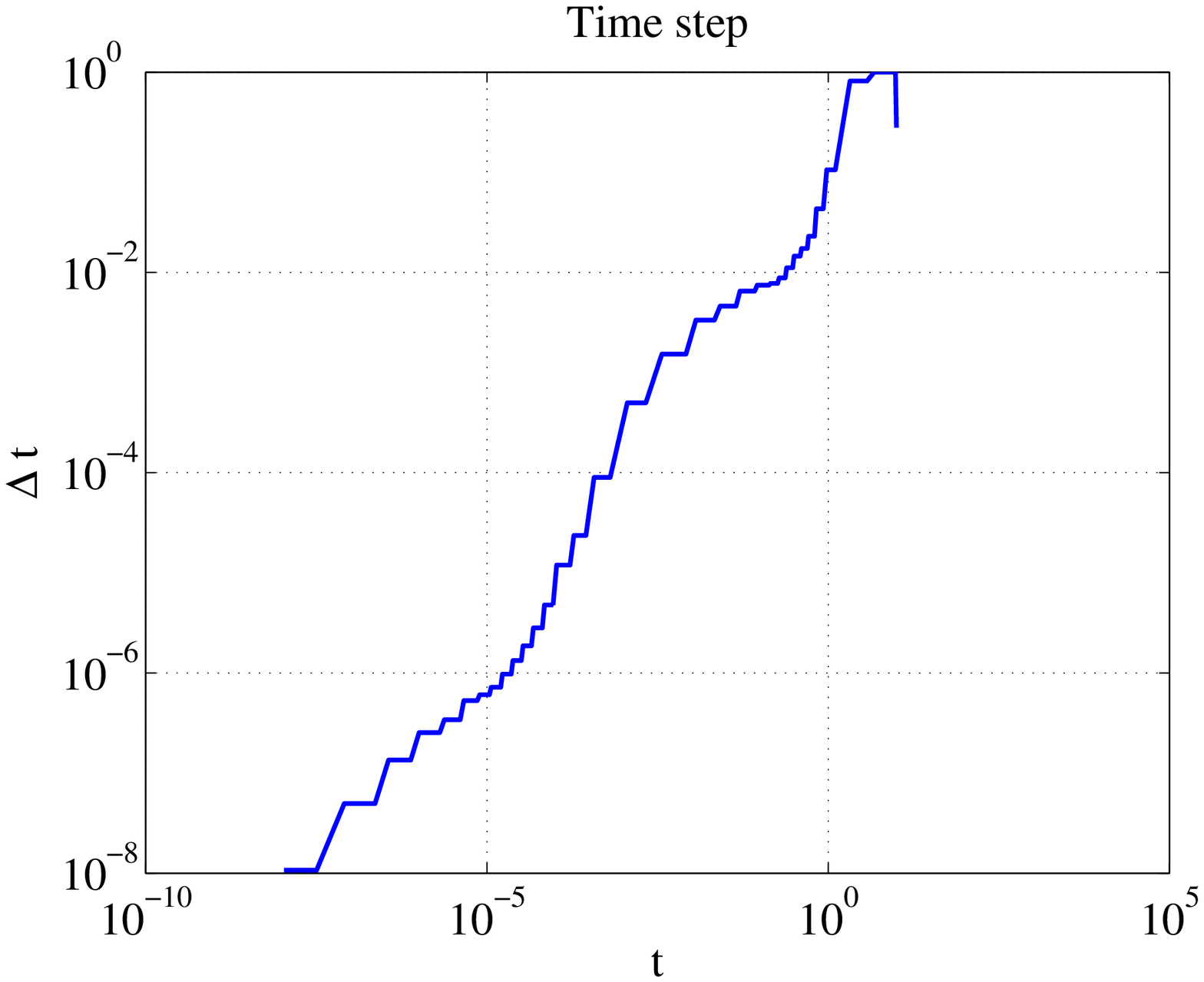} 
    \end{tabular}
    \caption{Mesh trajectories (top) and time step histories (bottom) for
      the given solution \en{ugiven1} (Example~1) with various values of
      $\tau$ using MMPDE6 and the arclength monitor function.}
    \label{fig1}
  \end{figure}

\end{example}

%% EXAMPLE 2.
\begin{example}
  A function which provides a more stringent test of a moving mesh
  calculation is  
  \begin{gather} 
    u(x,t) = e^{-\pi^2t} \sin(\pi x) + e^{-(10\pi)^2t} \sin(2\pi x),
    \label{eq:ugiven2}
  \end{gather}
  which is a slight modification of \en{ugiven1} having two
  widely-separated time scales over which the solution varies.  Notice
  in the results depicted in Figure~\ref{fig1b} that the more rapid
  variation embodied by the second term in \en{ugiven2} is only really
  well-captured in the mesh for the smallest value of $\tau=10^{-5}$.
  When the fast term dies out shortly after time $t=10^{-4} \;s$, then
  the mesh redistributes to resolve the remaining single peak.
  \begin{figure}[hbtp]
    \centering
    \psfragfontsize
    \psfrag{x}[c][b]{$x$}
    \psfrag{t}[Bc][c]{$t$}
    \psfrag{Mesh trajectories}{}
    \begin{tabular}{ccc}
      \figfont{$\tau = 10^{-1}$, $\text{CPU}=11\;s$} & 
      \figfont{$\tau = 10^{-3}$, $\text{CPU}=17\;s$} & 
      \figfont{$\tau = 10^{-5}$, $\text{CPU}=27\;s$} \\
      \includegraphics[width=0.3\textwidth]{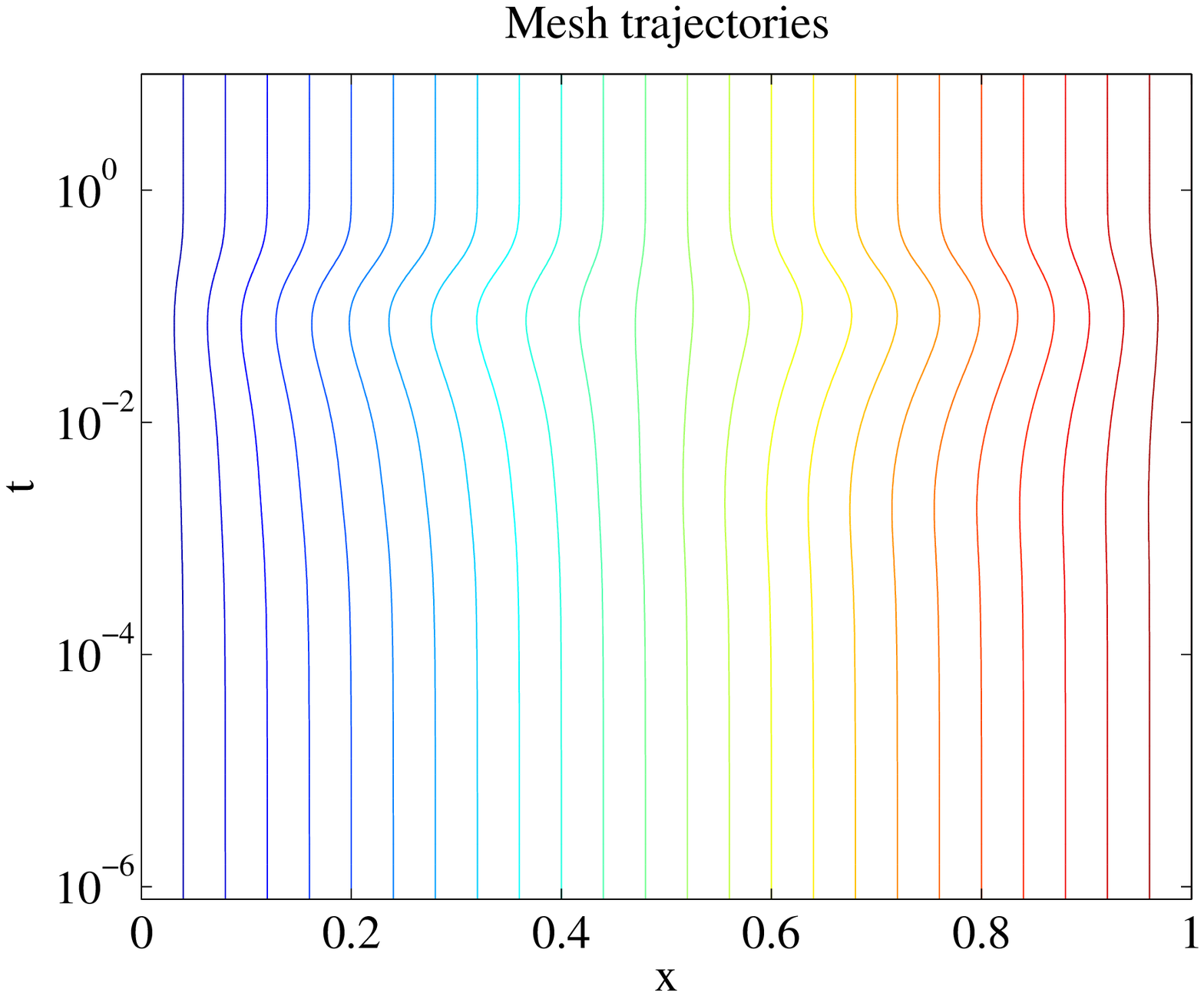} &
      \includegraphics[width=0.3\textwidth]{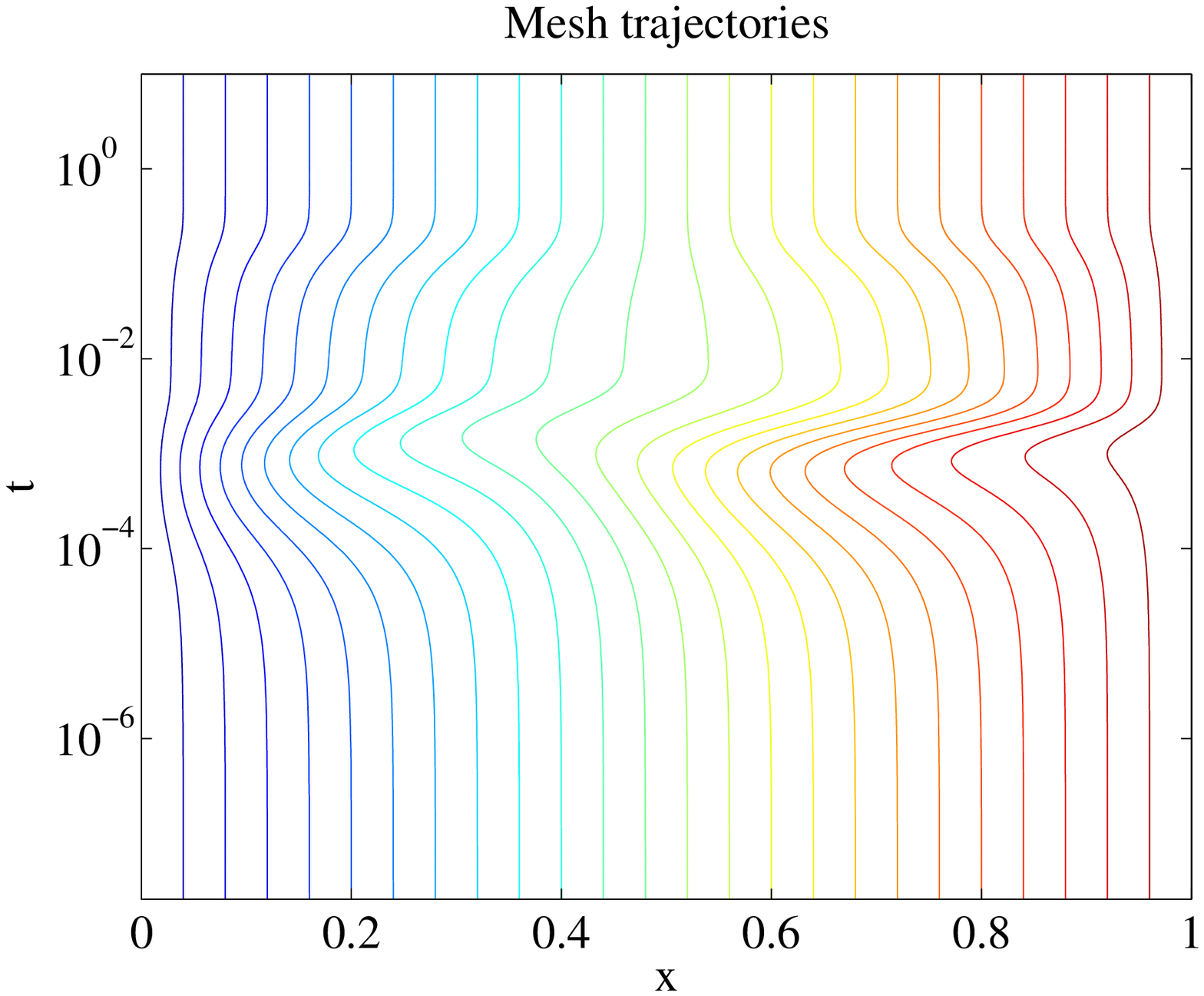} &
      \includegraphics[width=0.3\textwidth]{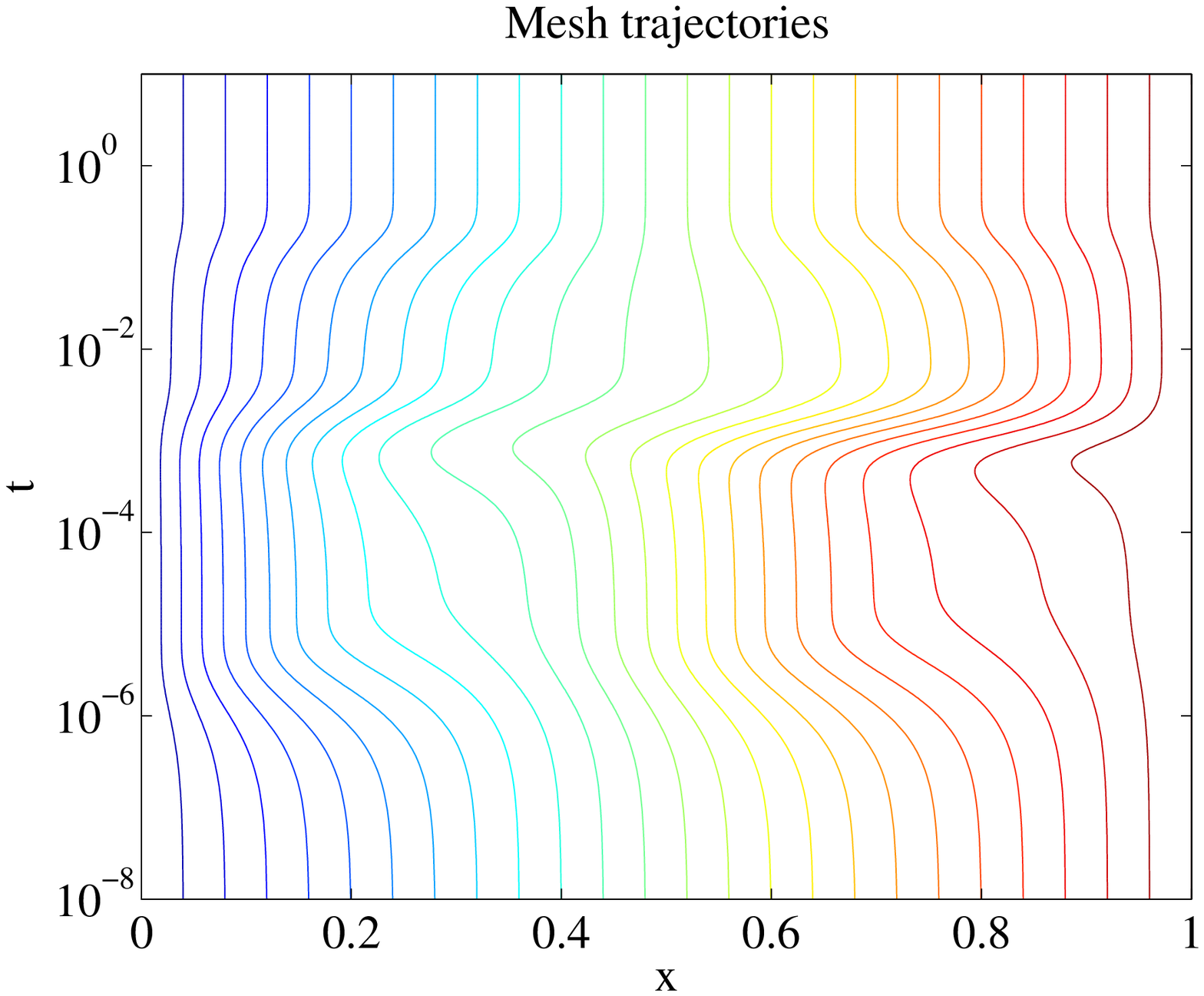}
%%      \\
%%      \psfrag{t}[c][b]{$t$}
%%      \includegraphics[width=0.3\textwidth]{u3-1-dt.eps} &
%%      \psfrag{t}[c][b]{$t$}
%%      \includegraphics[width=0.3\textwidth]{u3-3-dt.eps} &
%%      \psfrag{t}[c][b]{$t$}
%%      \includegraphics[width=0.3\textwidth]{u3-5-dt.eps} 
    \end{tabular}
    \caption{Mesh trajectories 
      %% (top) and time step histories (bottom)
      for the given solution \en{ugiven2} (Example~2) using MMPDE6 and
      the arclength monitor function.} 
    \label{fig1b}
  \end{figure}
\end{example}

%% EXAMPLE 3.
\begin{example}
  We next consider the following Gaussian function which blows up at the
  point $x=\xstar$ as $t\rightarrow \tstar$
  \begin{gather}
    u(x,t) = \frac{1}{\sqrt{4\pi(\tstar - t)}} \,
    \exp\left(-\frac{\beta(x-\xstar)^2}{4(\tstar-t)}\right) 
    \label{eq:exp}
  \end{gather}
  and which is more typical of solutions to the nonlinear diffusion
  equations we study later.  Here, we take $\beta = 100$ to ensure the
  blow-up region is very narrow and let $\xstar=0.5$ and $\tstar=0.4$.
  Typical solution curves are shown in Figure \ref{fig1u}, and the
  mesh trajectories and time step behaviour are displayed in 
  Figure~\ref{fig:exp} for a number of constant values of $\tau$ ranging
  between $10^{-1}$ to $10^{-5}$.  Note that the vertical axis for the
  mesh trajectories is displayed in terms of $(\tstar-t)$ on a log scale
  so that the clustering of mesh points near the blow-up time is
  actually visible.
  \begin{figure}[htbp]
    \centering
    \psfragfontsize
    \psfrag{x}[c][b]{$x$}
    \psfrag{t}{}
    \psfrag{ t}{}
    \psfrag{-t}[Bc][c]{$\tstar-t$}
    \psfrag{*}{}
    \psfrag{D}[Bc][c]{$\Delta t$}
    \psfrag{Mesh trajectories}{}
    \psfrag{Time step}{}
    \begin{tabular}{ccc}
      \figfont{$\tau=10^{-1}$} & \figfont{$\tau=10^{-3}$} &
      \figfont{$\tau=10^{-5}$} \\
      \includegraphics[width=0.3\textwidth]{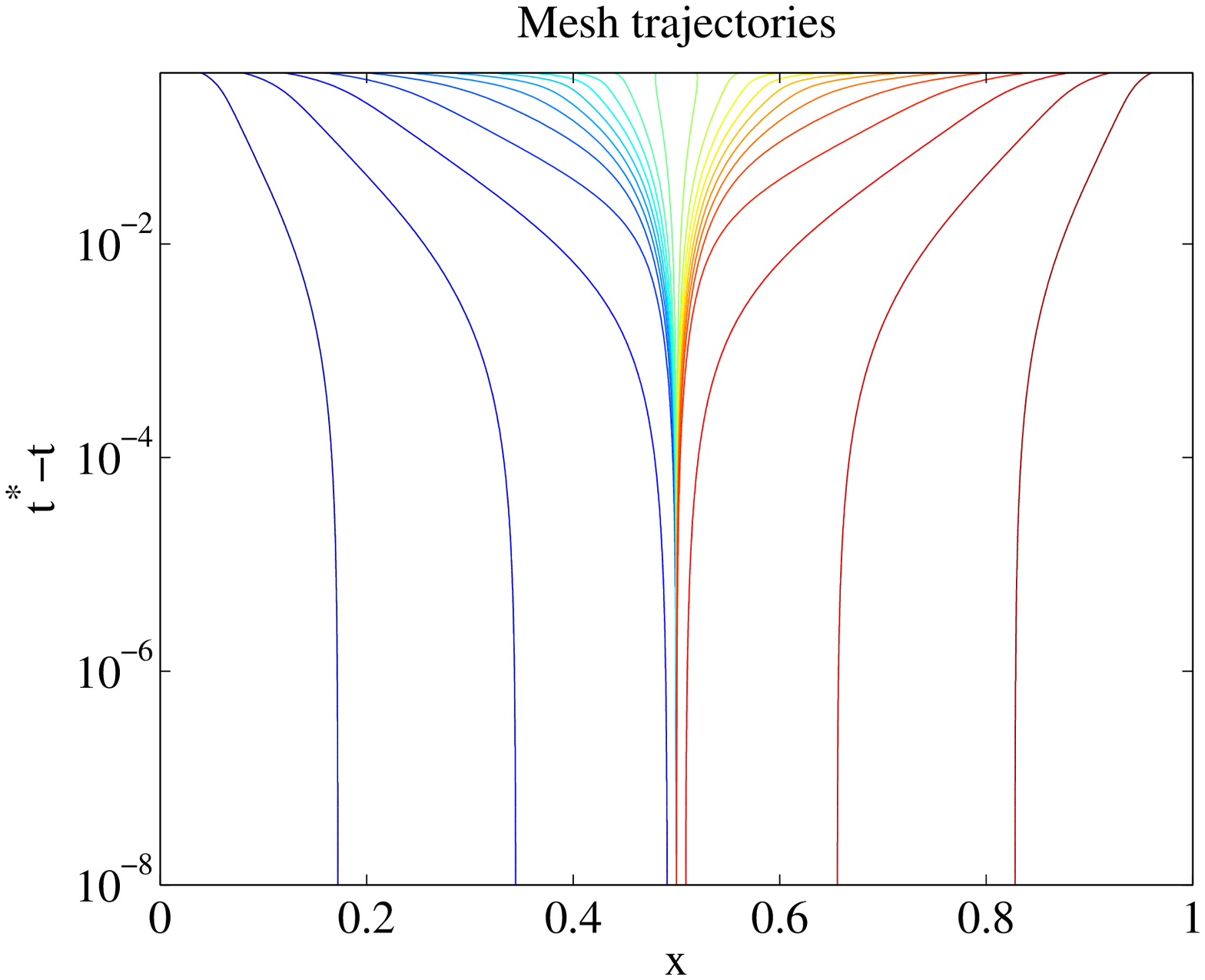} & 
      \includegraphics[width=0.3\textwidth]{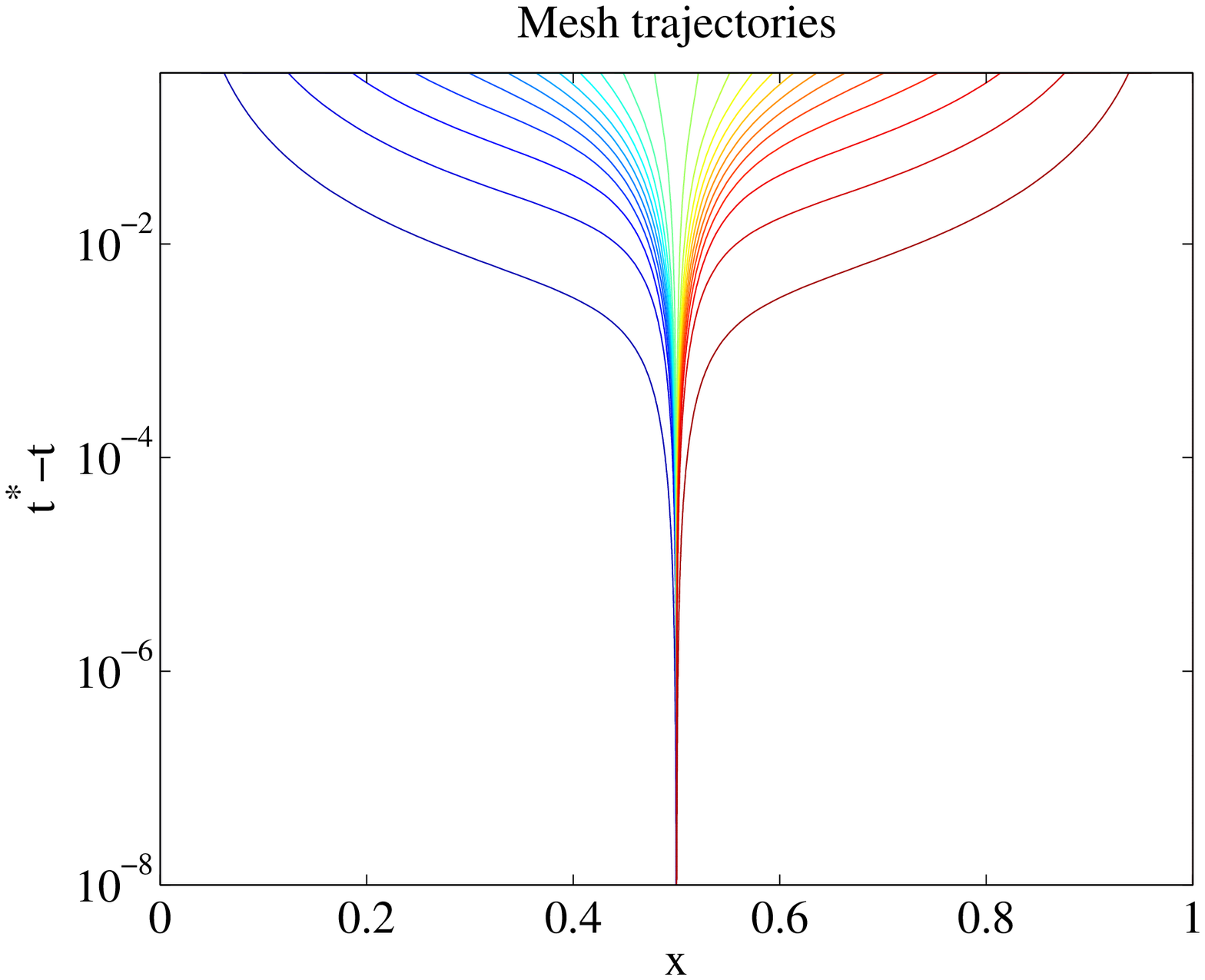} & 
      \includegraphics[width=0.3\textwidth]{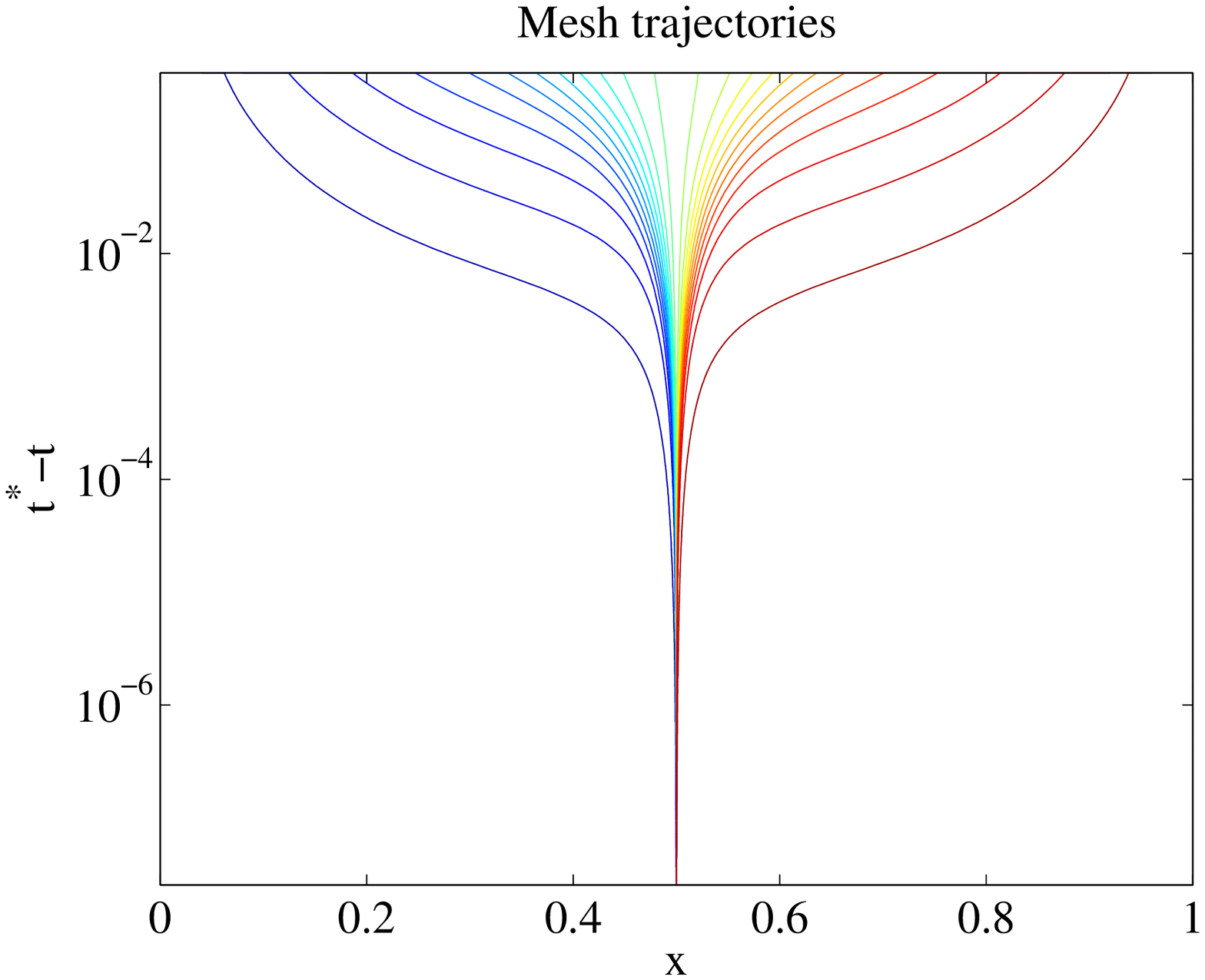} \\
      \psfrag{-t}[c][b]{$\tstar-t$}
      \includegraphics[width=0.3\textwidth]{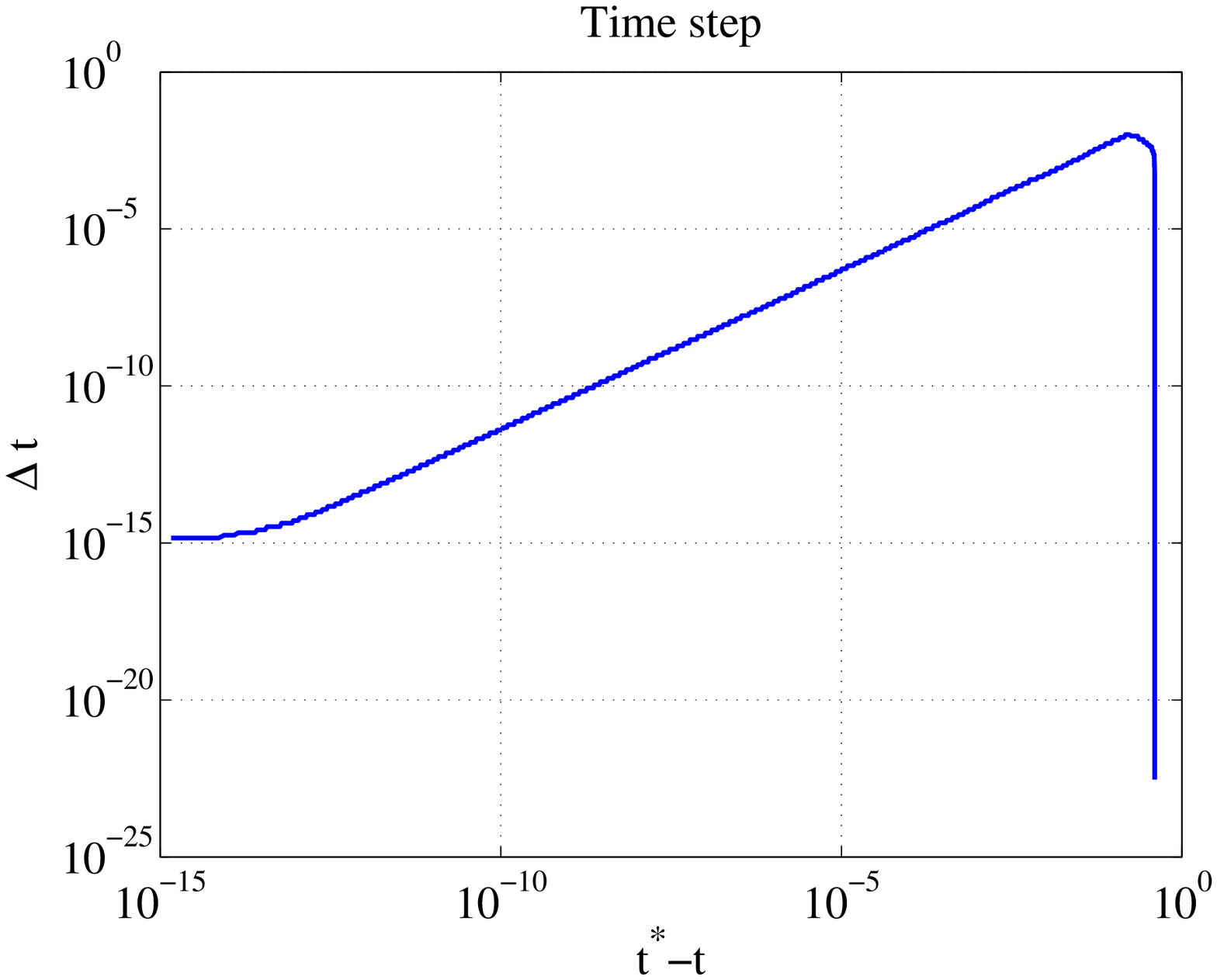} &
      \psfrag{-t}[c][b]{$\tstar-t$}
      \includegraphics[width=0.3\textwidth]{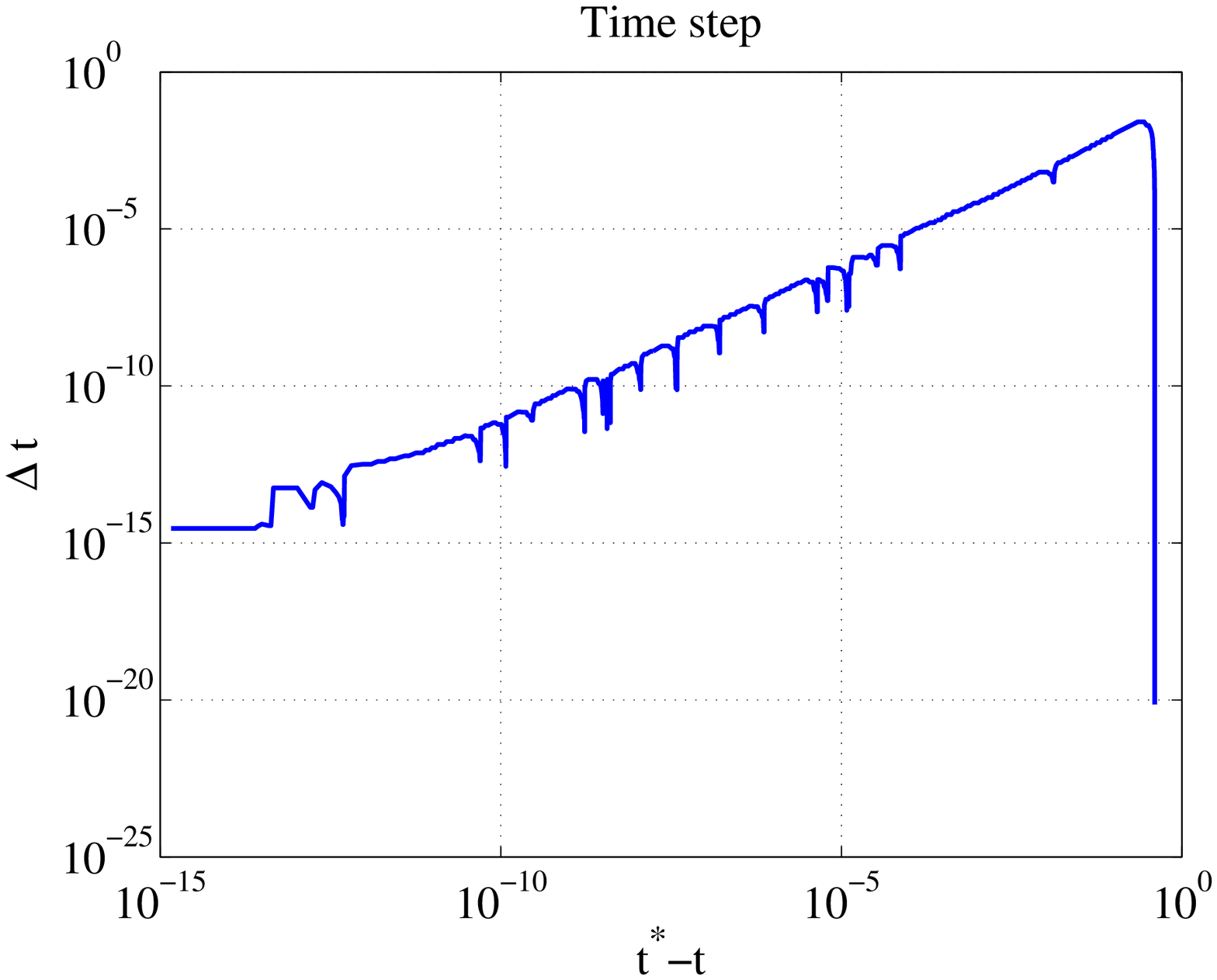} &
      \psfrag{-t}[c][b]{$\tstar-t$}
      \includegraphics[width=0.3\textwidth]{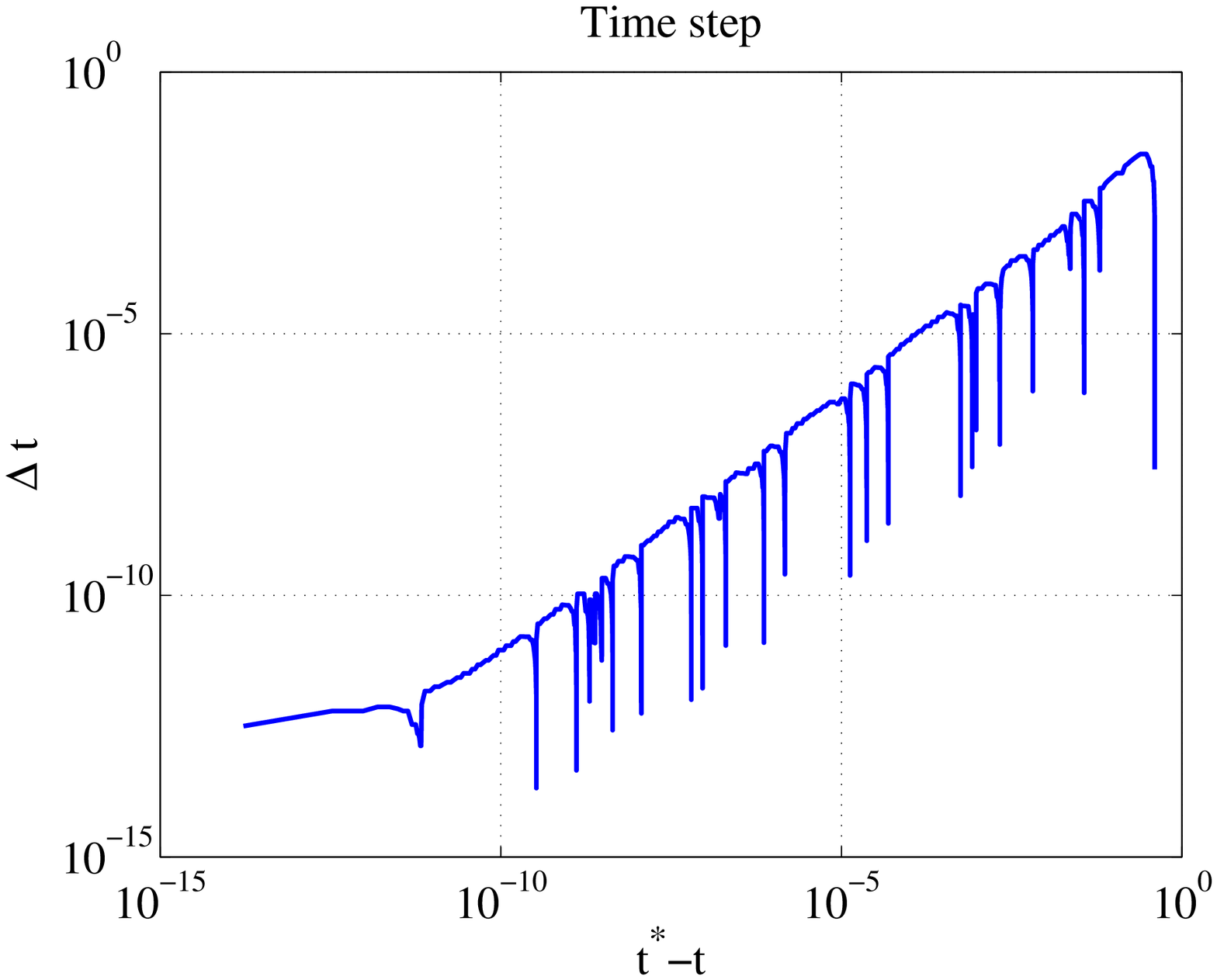} 
    \end{tabular}
    \caption{The mesh trajectories (top) and time step histories (bottom)
      for the Gaussian function \en{exp}.  Note that the horizontal axis
      measures $\tstar-t$, so that time progresses from right to left.}
    \label{fig:exp}
  \end{figure}
  When $\tau$ is taken as large as $10^{-1}$, there is clearly
  sufficient smoothing in the mesh that a significant number of mesh
  points transfer outside the blow-up region.  As $\tau$ is reduced in
  size, the mesh is closer to being equidistributed and the resolution
  of the blow-up peak is much sharper. However, this higher
  concentration of mesh points comes at the expense of a stiffer mesh
  equation, as evidenced by a reductions in the allowable time step in
  DDASSL and numerous time step failures.  As before, once $\tau$ is
  taken smaller than $10^{-5}$, there is no longer any visible
  difference in the computed solution.
\end{example}

\leavethisout{
  \begin{example}
    As a final example, we consider a function with a sharp front that
    propagates in time which was proposed in~\cite{hrr94b}:
    \begin{gather*}
      u(x,t) = \frac{1}{2} \left[ 1 - \tanh \left(
          c(t) (x-t-0.4) \right) \right],  
    \end{gather*}
    where 
    \begin{gather*}
      c(t) = 1 + \frac{1}{2}(10^3-1)\, \left[ 1 + \tanh \left(
          100 (t-0.2) \right) \right],  
    \end{gather*}
    for $0\leq x \.eq 1$ and $0 \leq t \leq 0.55$.  The solution forms a
    sharp front suddenly at about time $t=0.2$, which then propagates to
    the right at speed 1.
  \end{example}
}

\paragraph*{Remark.}
The limitations on $\tau$ indicated in the above three examples are
problem-dependent.  Until the present time, all moving mesh calculations
appearing in the literature have been performed with a constant value of
$\tau$ which in practice must essentially be chosen by trial and error.
Because $\tau$ represents a time scale for the mesh evolution, it is not
appropriate to take $\tau$ constant when the solution undergoes rapid
changes that require the mesh to respond on very different time scales,
such as might occur in the case of blow-up or shock motion with highly
variable front speeds.  In these situations, it makes much more sense to
vary $\tau$ over time in a way that adapts the mesh throughout a
computation so as to respond over a suitable time scale to changes in
the solution.  We will examine the issue of choosing an appropriate form
of $\tau(t)$ in the context of blow-up problems, which are described
further in the next section.

\leavethisout{
  \subsection{Mesh quality measures}
  
  There are a number of mesh quality measures which are introduced in
  \cite{hrr94b} and which could be employed here:
  \begin{itemize}
  \item $\ds{E(t) :=  \max_{2\leq i \leq N-1} |E_i(t)|}$
    where $E_i$ are defined in \en{discrete-rhs}, which is a measure of the
    deviation of the grid points from an exactly equidistributed mesh.
  \item $\ds{G(t) :=  \sum_{i=0}^{N-1} \frac{1}{2} \left((M_t)_{i+1} +
        (M_t)_{i}\right)\, (x_{i+1}-x_i) \approx \int_0^1
      \frac{\partial}{\partial t} M(x,t)\, dx}$.
  \item $\ds{S(t) :=  \max_{0\leq i \leq N} |\dot{x}_i(t)|}$.
  \item $\ds{L(t) :=  \max_{0\leq \xi \leq 1}
      \frac{M(x(\xi,0),0)}{M(x(\xi,t),t)}},$ which is a measure of the
    stability of the mesh.  
  \end{itemize}
}
  
%%%%%%%%%%%%%%%%%%%%%%%%%%%%%%%%%%%%%%%%%%%%%%%%%%%%%%%%%%%%%%%%%%%%%%%%
\section{Self-similar blow-up}
\label{sec:blowup}

An ideal class of problems with which to examine the behaviour of moving
mesh methods is that which models blow-up phenomena.  One of the
simplest equations in this class, and one which will form the basis of
most of the numerical simulations presented in this paper, is the
following nonlinear diffusion equation of parabolic type
\cite{bcr05,bhr96}:
\begin{subequations}\label{eq:blowup}
  \begin{gather}
    u_t=u_{xx} + u^p, \label{eq:blowup1}\\
    \intertext{with boundary and initial conditions}
    u(0, t)= u(1, t)=0 \quad \text{and} \quad
    u(x,0)=u_0(x). \label{eq:blowup2}  
  \end{gather}
\end{subequations}
This equation models, for example, the temperature in a reacting medium.
It is well-known \cite{bb92,fm85} that if $u_0(x)$ is sufficiently
large, positive, and has a single non-degenerate maximum, then there is
a blow-up time $\tstar < \infty $ and a unique blow-up point $\xstar$
such that
\begin{gather*} 
  u(\xstar,t)\longrightarrow \infty \quad {\rm as}\quad
  t\longrightarrow \tstar, \\
  \intertext{and}
  u(x,t)\longrightarrow u(x,\tstar)<\infty
  \quad {\rm if} \quad x \neq \xstar.
\end{gather*}
That is, even when there are smooth initial data the solution becomes
unbounded at an isolated point $\xstar$ in finite time.  Other forms of
the nonlinear term in \en{blowup1} will also lead to blow-up (for
example, with the nonlinear term $u^p$ replaced by $e^u$) but the
polynomial form is particularly convenient 
for our purposes because of its scaling properties, which we describe
next.

This equation has been very well-studied in the mathematical literature
and the solutions are known to exhibit self-similar behaviour.  In
particular, if we define $\beta=1/(p-1)$, then the solution has a
self-similar profile which blows up according to
\begin{gather}
  u \sim (\tstar - t)^{-\beta} \label{eq:uasympt}
\end{gather}
asymptotically as $t\longrightarrow \tstar$~\cite{bb92}.  This
information was used by Budd\ \etal~\cite{bhr96} as part of a scaling
argument 
to show that the MMPDE corresponding to the blow-up problem \en{blowup}
is also scale-invariant if the monitor function is
chosen to be $M=u^{p-1}$.  They then presented a series of numerical
simulations which showed that the MMPDE method is capable of reproducing
the self-similar solution profiles in a more accurate and stable manner
than is possible with other more common choices of monitor function such
as arclength.  

An essential observation made in \cite{bhr96}, which has particular
importance for this paper, is that the mesh in the MMPDE method has a
natural time scale that is determined by scaling arguments.  If the
scale-invariant monitor function $M=u^{p-1}$ is employed in
calculations, then we know from \en{uasympt} that $M \sim
(\tstar-t)^{-1}$ asymptotically as $t\longrightarrow \tstar$.  It is
then straightforward to show that the mesh has a natural time scale of
motion which is determined by the choice of the MMPDE; in particular,
\newcommand{\Tmesh}{\textsub{T}{mesh}}
\begin{align*}
  \Tmesh &= \order{\tau} & & \text{(for MMPDE4),}\\
  \intertext{and}
  \Tmesh &= \order{\frac{\tau}{M}} \sim \tau(\tstar - t)
  & & \text{(for MMPDE6).}
\end{align*} 
Budd \etal\ argue in the first case that when $\tau$ is taken to be a
constant, the ability of the mesh to react to changes in the solution is
limited by the lower bound $\tau$ on the mesh time scale and so MMPDE4
does not allow the mesh to evolve all the way into the blow-up.  In
other words, once $|\tstar - t| < \tau$, the mesh will no longer evolve
rapidly enough to keep up with the solution.  On the other hand, MMPDE6
does allow the mesh to evolve even when $t$ is close to $\tstar$,
because of the extra factor of $(\tstar-t)$ appearing in $\Tmesh$.
This hypothesis regarding the superiority of MMPDE6 over MMPDE4 for the
blow-up problem \en{blowup} is borne out in computations~\cite{bhr96}
where MMPDE6 is capable of capturing the self-similar solution profile
much further into blow-up than MMPDE4.

\subsection{A strategy for varying $\tau$}
\label{sec:varytau}

Our main claim in this paper is that requiring a small,
constant value of $\tau$ can introduce unnecessary stiffness in the
moving mesh PDE.  In the case of solutions to \en{blowup}, blow-up
occurs at a point $\xstar$ which is stationary, and so there is an
initial transient mesh motion in which mesh points race into the blow-up
region, after which the mesh points are relatively stationary even
though the solution $u$ and the monitor function $M$ are both increasing
rapidly.  It is therefore natural to suggest that capturing the initial
mesh transients may require a small initial value of $\tau$, but that
$\tau$ can be significantly increased later on in the blow-up process at
little risk of negatively impacting the accuracy of the mesh locations.
Recalling the examples considered in Section~\ref{sec:tau-ex}, we
reiterate that decreasing $\tau$ allows the mesh to react to solution
changes more rapidly, but also introduces additional stiffness into the
MMPDE; conversely, increasing $\tau$ speeds up the computations but may
unnecessarily smooth out the mesh and adversely affect solution
accuracy.  Adapting $\tau$ as described above should therefore act to
minimize the stiffness in the MMPDE in later stages of blow-up and so
decrease computational cost.

With this in mind, we propose the following solution-adaptive strategy
for choosing $\tau$:
\begin{itemize}
\item Set $\tilde{\tau}(t) = \tau_o \,\ds{\max_i} (M_i)$,\ 
  where $\tau_o$ is a constant.
  %%\label{eq:tau-mon}
\item Choose $\tau(t) = \min \left( \max \left[ \tilde{\tau}(t), 
      \textsub{\tau}{min} \right], \textsub{\tau}{max} \right)$,  
  which forces $\tau$ to lie in the interval $[\textsub{\tau}{min}, 
  \textsub{\tau}{max}]$. 
\end{itemize}
This ensures that the mesh time scale $\Tmesh$ is small in the initial
stages of blow-up, but increases to $\textsub{\tau}{max}$ later on when
the mesh velocities are much smaller.  It is important to point out that
this strategy is applicable \emph{only} to blow-up problems with
self-similar structure of this sort, and not for more general
situations.   

\leavethisout{
  \subsection{Alternate choices for $\tau$}
  
  For problems that don't have the special structure of self-similar
  blow-up problems, a more general strategy for selecting $\tau$ is
  required.  The main principle in selecting  the mesh relaxation time
  is that $\tau$ should be small (and hence introduce less 
  smoothing) when the mesh moves very rapidly so as to keep the mesh much
  closer to equidistributed.  On the other hand, $\tau$ should be allowed
  to take on larger values when the mesh is not evolving as quickly.
  Consequently, the we make the following heuristic assumption, that
  $\tau$ vary depend on the mesh velocity according to 
  \begin{gather}
    \tau = \frac{\alpha}{\max_x |\dot{x}|},
    \label{eq:tau-xdot}
  \end{gather}
  where $\alpha$ is some constant.  The mesh velocity can be written in
  terms of easily computed quantities as
  \begin{gather*}
    \dot{x} = \frac{M_t}{M_\xi/x_\xi},
  \end{gather*}
  which can then be approximated at discrete points by
  \begin{gather}
    \dot{x}_i \approx
    \frac{(M_i^{n+1} - M_i^n)\, (x_{i+1}^n-x_{i-1}^n)}{\Delta t\,
      (M_{i+1}^n-M_{i-1}^n)}. 
    \label{eq:tau-xdot-discrete}
  \end{gather}
}

There are a number of other approaches for selecting an appropriate mesh
time scale which have been developed in the context of other moving mesh
methods (see \cite{af86b,hl86}).  However, we have found that neither of
these approaches is effective for the blow-up problems under
consideration here. 

\leavethisout{
  \begin{mynote}
    Some other choices:
    \begin{itemize}
    \item Hyman and Larrouturou~\cite{hl86} suggested a mesh time scale
      which is based on the time variation of the solution:
      \begin{gather}
        \tau^{-1} = \max \left( \textsub{\tau}{max}^{-1}, 
          \ds{\min_i} \left[ \frac{F(U_i)}{\ds{\max}( |U_i|, \textsub{U}{min})}
          \right] \right). 
      \end{gather}
      Here, both $\textsub{\tau}{max}$ and $\textsub{U}{min}$ are positive
      constants that ensure $\tau^{-1}$ is well-defined when the solution
      reaches steady state or the solution components vanish.  The right hand
      side functions $F$ are approximated with finite differences in $U$.
    \item Adjerid and Flaherty~\cite{af86b} propose another choice of mesh
      time scale.
    \item Also describe the approach in \cite{huang01}.
    \end{itemize}
  \end{mynote}
}

%%%%%%%%%%%%%%%%%%%%%%%%%%%%%%%%%%%%%%%%%%%%%%%%%%%%%%%%%%%%%%%%%%%%%%%%
\section{Numerical experiments}
\label{sec:results}

We now consider a number of computational examples in which the mesh
equation is coupled with the physical PDE.  When the blow-up problem
\en{blowup} is transformed into a moving coordinate system, it can be
written in the following form
\begin{gather*}
  \dot{u} - \frac{u_\xi}{x_\xi} \, \dot{x} = \frac{1}{x_\xi} \left(
  \frac{u_\xi}{x_\xi} \right)_\xi + u^p. 
\end{gather*}
We employ a method-of-lines approach in which this equation is
discretized with second order spatial accuracy using centered finite
differences to obtain the following equation for the solution values
$u_i(t)$: 
\begin{gather}
  \dot{u}_i - \frac{u_{i+1}-u_{i-1}}{x_{i+1}-x_{i-1}} \,\dot{x}_i = 
  \frac{2}{x_{i+1}-x_{i-1}} \left(
    \frac{u_{i+1}-u_i}{x_{i+1}-x_i} - 
    \frac{u_{i}-u_{i-1}}{x_{i}-x_{i-1}} \right) + u_i^p,
  \label{eq:discrete-pde}
\end{gather}
where the ``dot'' refers to a time derivative.  The resulting coupled
system of nonlinear ODEs which governs the mesh and solution,
\en{discrete-mmpde} and \en{discrete-pde}, is then integrated in time
using the stiff ODE solver DDASSL~\cite{pet82} with a finite difference
Jacobian.  Unless indicated otherwise, we use absolute and relative
error tolerances of $10^{-8}$.  Homogeneous boundary conditions are
imposed so that $u_0(t)= u_N(t)=0$, and we take the initial solution
profile $u(x,0)=20 \sin(\pi x)$.  The initial mesh, $x_i(0)$ is
determined by equidistributing based on the initial conditions.
Computations are performed using blow-up exponents $p=2$ and $p=5$, and
the monitor function is taken to be $M=|u|^{p-1}$, which preserves
scaling invariance of the mesh.  The variable $\tau$ simulations are
performed using $\tau(t) = 10^{-8}\,\max_x(M)$ (that is,
$\tau_o=10^{-8}$) and then enforcing that $\tau$ lie in the interval
$[10^{-8}, 10^{-1}]$.

One aim of these computations is to compute as far into blow-up as
possible and to obtain the best possible estimate of the blow-up time
$\tstar$.  In all our simulations, we compute as far as DDASSL will
allow, up until such time as the solver fails (which in practice
manifests itself as a time step selection failure).

Since no exact analytical solution is available for this problem, it is
difficult to assess the accuracy of a given computed solution.  In this
paper, we employ a number of qualitative and quantitative measures to
compare the accuracy of the computed solutions:
\begin{itemize}
\item The termination time (as an estimate of $\tstar$) is compared to
  the blow-up time determined from a highly-resolved calculation, which
  gives a combined measure of the accuracy of the solution and the mesh.
  In particular, our best estimates of the blow-up times are
  $\tstar\approx 0.08243786$\ 
%%$\tstar\approx 0.082437318$\ 
  for $p=2$ and $\tstar\approx 1.5625962 \times 10^{-6}$\ 
%%  for $p=2$ and $\tstar\approx 1.5625950 \times 10^{-6}$\ 
%%for $p=2$ and $\tstar\approx 1.56259486 \times 10^{-6}$\ 
  for $p=5$, both of which are calculated with $N=1000$ points,
  variable $\tau$ with $\tau_o=10^{-8}$, and DDASSL error tolerances of
  $10^{-10}$. 
\item The value of $\ds{\textsub{u}{max}:=  \ds{\max_x}\, u}$, which is an
  indirect measure of solution accuracy, that represents how far the
  code is capable of computing into the singularity.  Ideally, we aim
  for $\textsub{u}{max}$ to be as large as possible.
\item The self-similarity of the various solution profiles computed over
  time is most easily determined by comparing directly to the following
  asymptotic formula derived in \cite{bhr96}:
  \begin{gather}
    \left(\frac{u}{\textsub{u}{max}}\right)^{p-1} \sim 
    \cos^2 \left( \pi(\xi - \textstyle{\frac{1}{2}}) \right) .
    \label{eq:ucos}
  \end{gather}
\end{itemize}
The computational cost of all subsequent simulations is compared by
measuring the elapsed CPU time on a 3 GHz Intel Xeon machine.

\subsection{Blow-up with $p=2$}
\label{sec:results-p2}

\begin{figure}[htbp]
  \centering
  \psfragfontsize
  \psfrag{xi}[c][b]{$\xi$}
  \psfrag{u/umax}[Bc][c]{$u/\textsub{u}{max}$}
  \psfrag{tstar-t}[c][b]{$\tstar-t$}
  \psfrag{x}[Bc][c]{$x$}
  \psfrag{Mesh trajectories}{}
  \includegraphics[width=0.4\textwidth,clip]{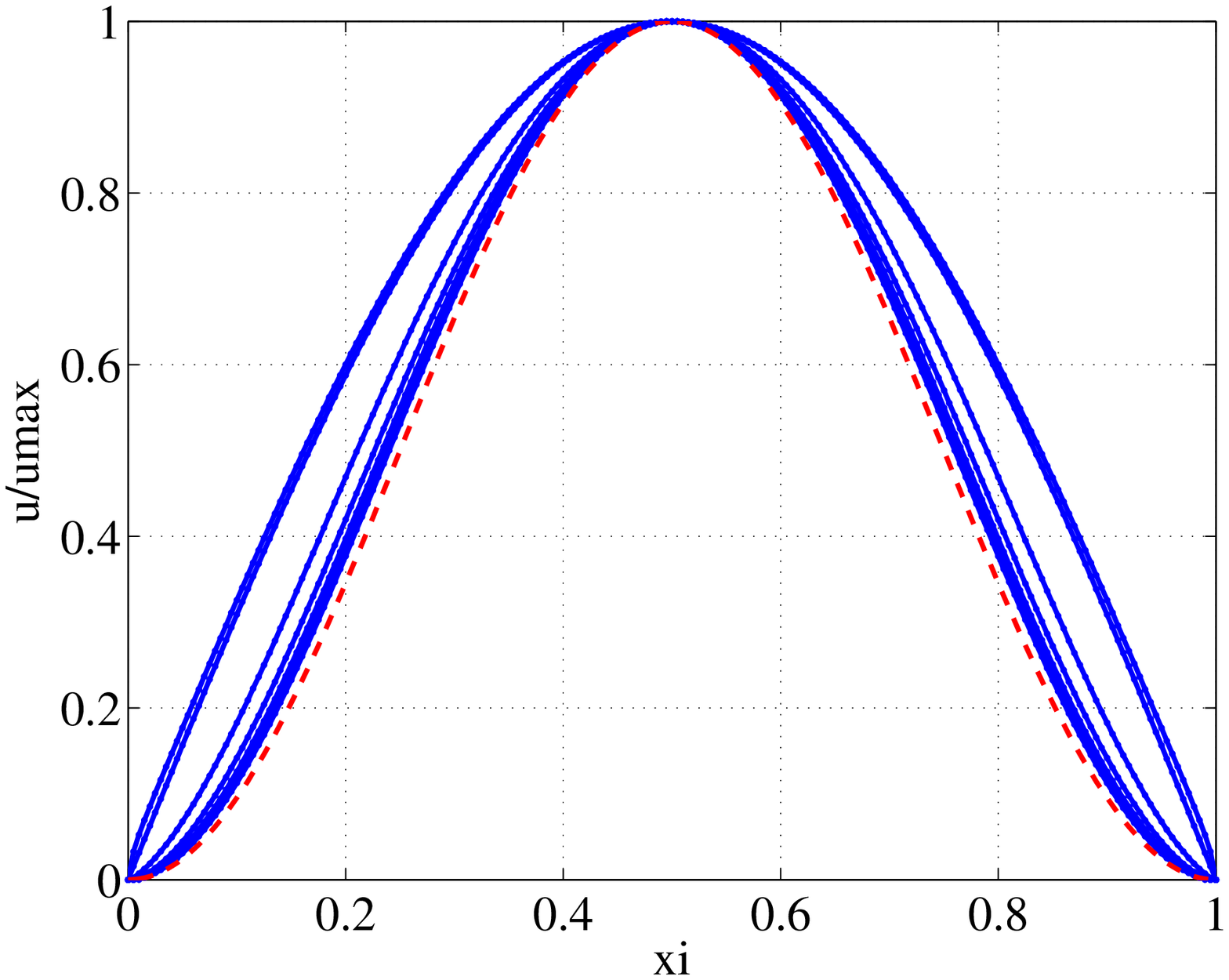} 
  \quad
  \includegraphics[width=0.4\textwidth,clip]{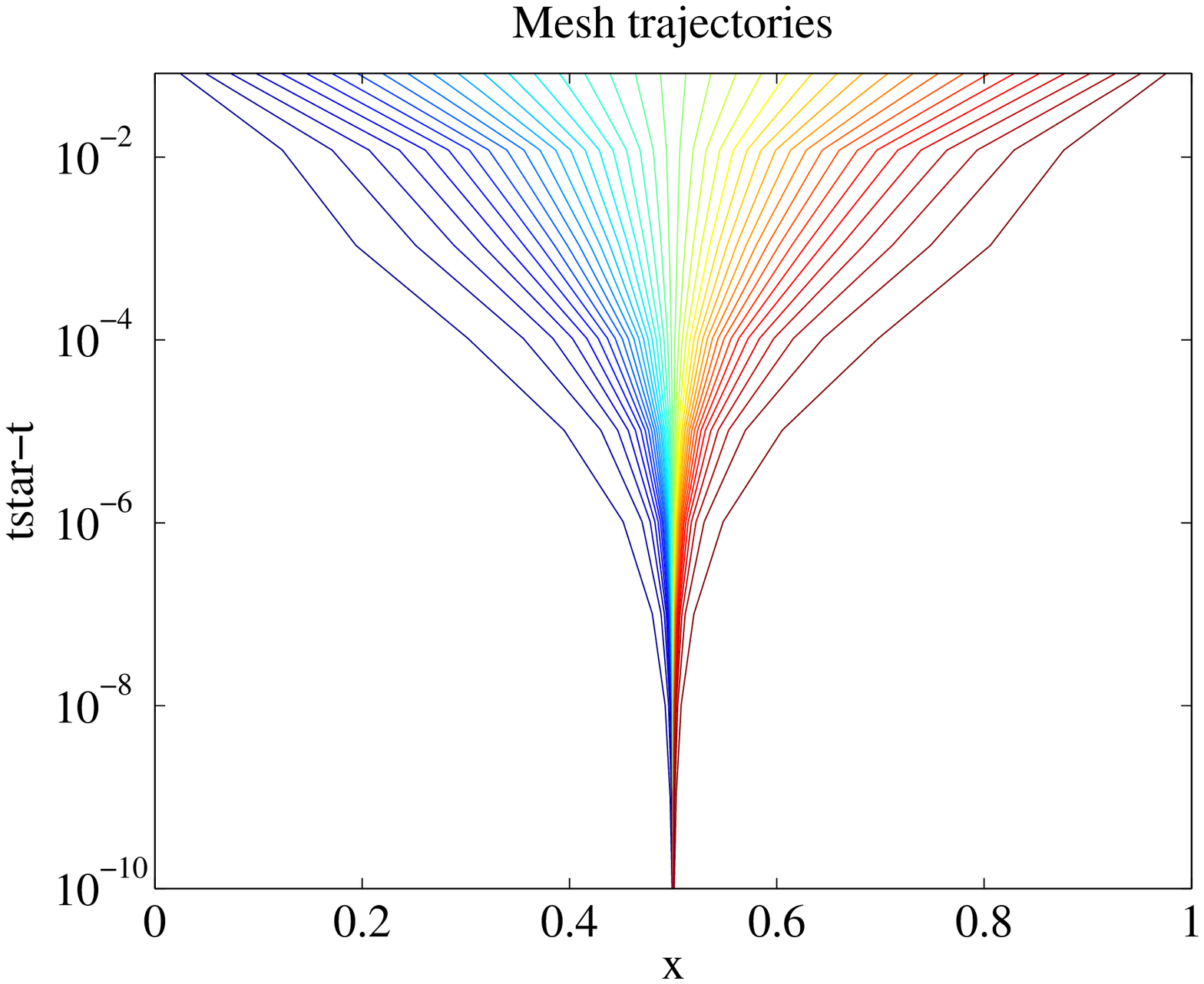} 
  \caption{Plot of the solution profiles (left) and mesh contours (right)
    for the blow-up problem with $N=200$ and $\tau=10^{-5}$ in the case
    $p=2$.  The self-similar profile is displayed as a dashed line for
    comparison.}
  \label{fig:p2const}
\end{figure}

We begin first by simulating the blow-up problem \en{blowup} with $p=2$.
For all computations in this section, we have not performed any mesh
smoothing (i.e., $ip=0$) in order to ensure that the computed solution
and mesh are as close to the similarity solution as possible.  The
solution and mesh contour plots for $N=200$ points are displayed for a
constant value of $\tau=10^{-5}$ in Figure~\ref{fig:p2const}.  The
various solution profiles correspond to a sequence of snapshots at times
where $\textsub{u}{max} = 10^n$, $n=1,2,\dots, 13$.  The plot of
$u/\textsub{u}{max}$ demonstrates how MMPDE6 with the monitor $M=u$ is
capable of capturing the self-similar nature of the solution.

\begin{figure}[htbp]
  \centering
  \psfragfontsize
  \psfrag{xi}[c][b]{$\xi$}
  \psfrag{u/umax}[Bc][c]{$u/\textsub{u}{max}$}
  \psfrag{tstar-t}[c][b]{$\tstar-t$}
  \psfrag{x}[Bc][c]{$x$}
  \psfrag{Mesh trajectories}{}
  \includegraphics[width=0.4\textwidth,clip]{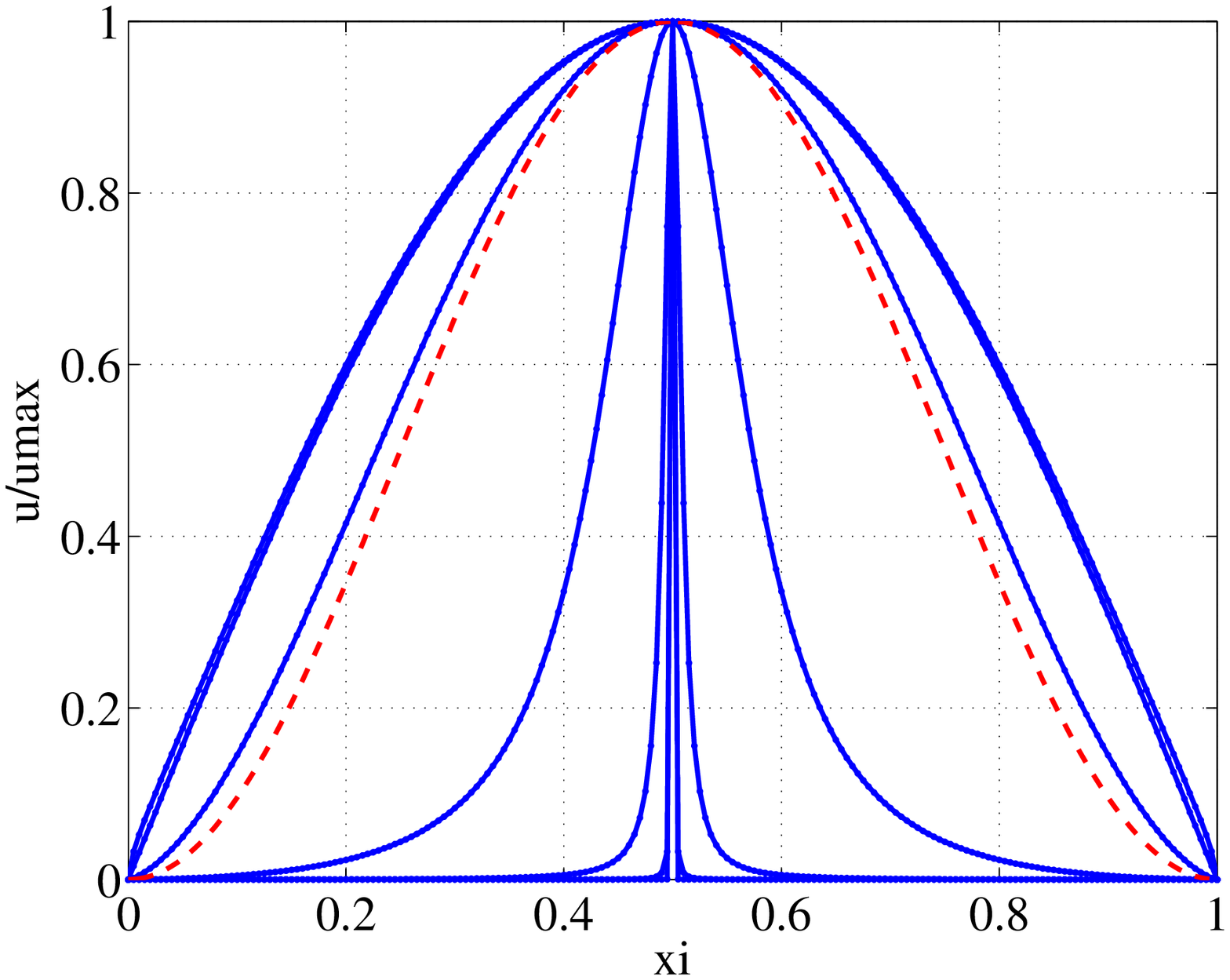}
  \quad
  \includegraphics[width=0.4\textwidth,clip]{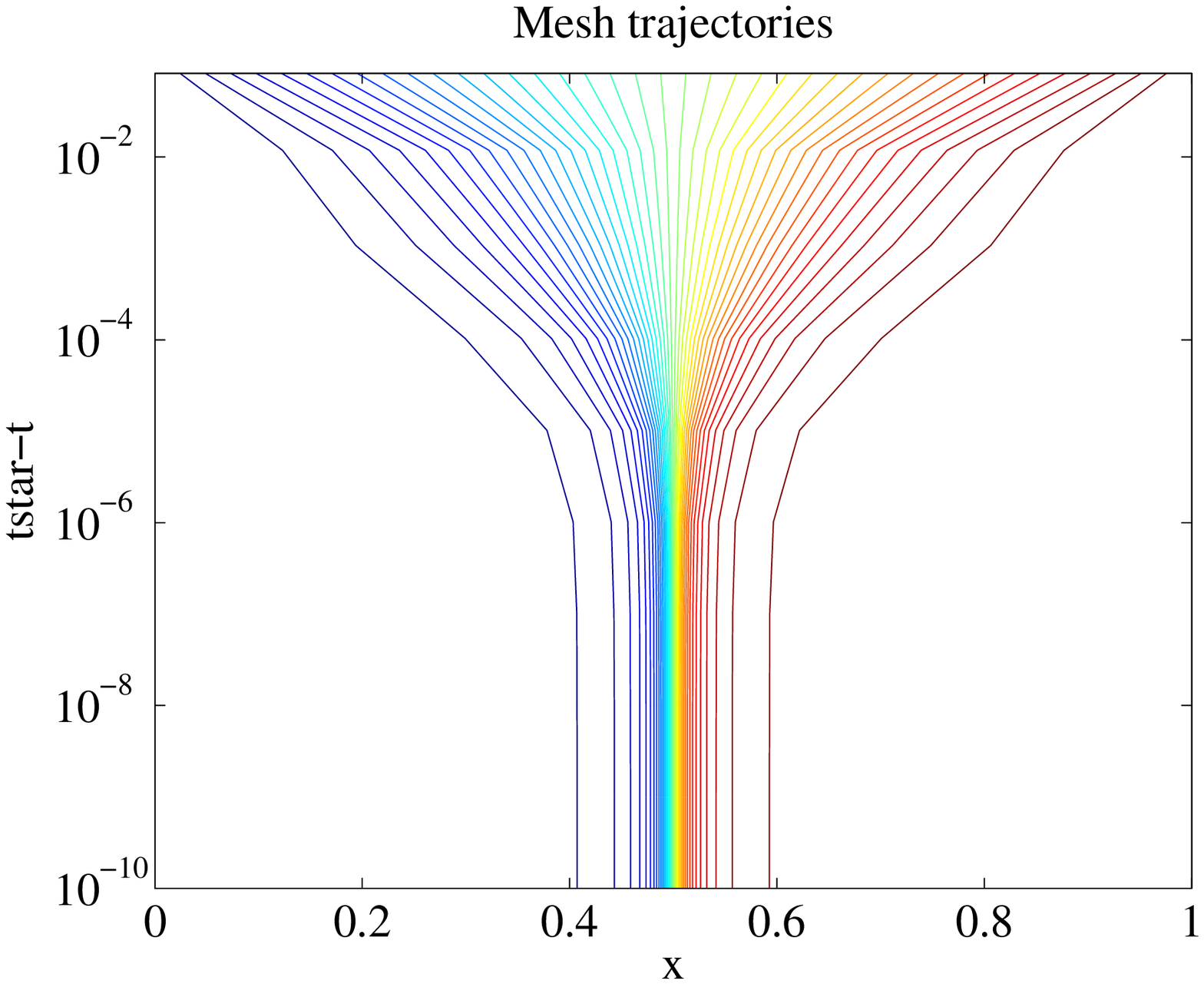}
  \caption{Plot of the solution profiles (left) and mesh contours (right)
    using MMPDE4, for the same parameters as in {\protect
      Figure~\ref{fig:p2const}}.} 
  \label{fig:p2const-mmpde4}
\end{figure}

To illustrate the effect of the choice of MMPDE on the solution, we have
also displayed the results for the same input data using MMPDE4 in
Figure~\ref{fig:p2const-mmpde4}.  This computation is clearly incapable
of maintaining grid resolution within the blow-up peak; in fact, by the
end of the calculation, the mesh degenerates to the extent that there
remains only a single grid point left to resolve the peak.  Furthermore,
this simulation fails at time $t=0.08243526\;s$, which is a much less
accurate estimate of the blow-up time than in the MMPDE6 calculations,
as we will see shortly.  Consequently, MMPDE6 is employed in the
remainder of the simulations in this paper.

\begin{figure}[htbp]
  \centering
  \psfragfontsize
  \psfrag{xi}[c][b]{$\xi$}
  \psfrag{u/umax}[Bc][c]{$u/\textsub{u}{max}$}
  \psfrag{tstar-t}[c][b]{$\tstar-t$}
  \psfrag{x}[Bc][c]{$x$}
  \psfrag{Mesh trajectories}{}
  \includegraphics[width=0.4\textwidth,clip]{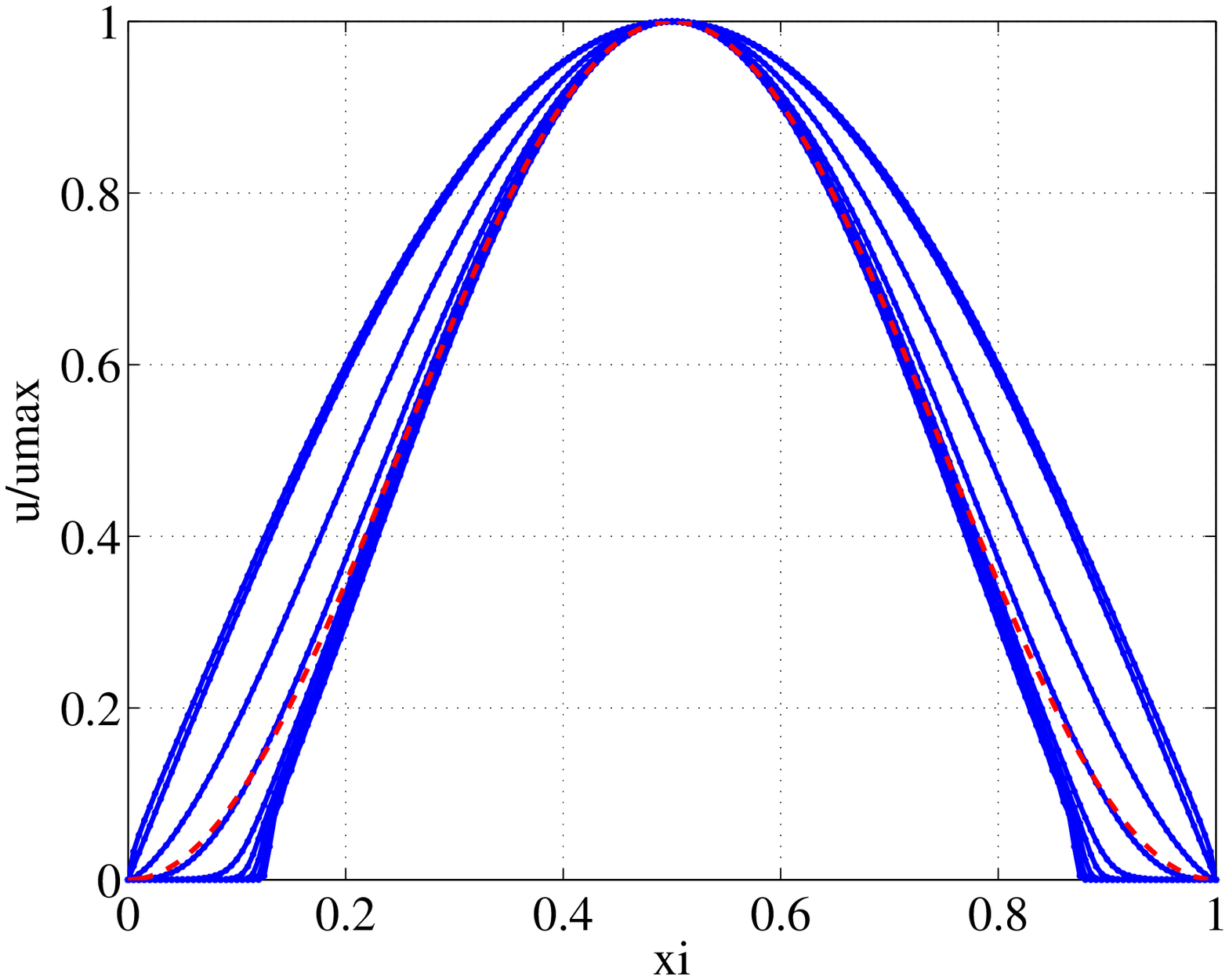}
  \quad
  \includegraphics[width=0.4\textwidth,clip]{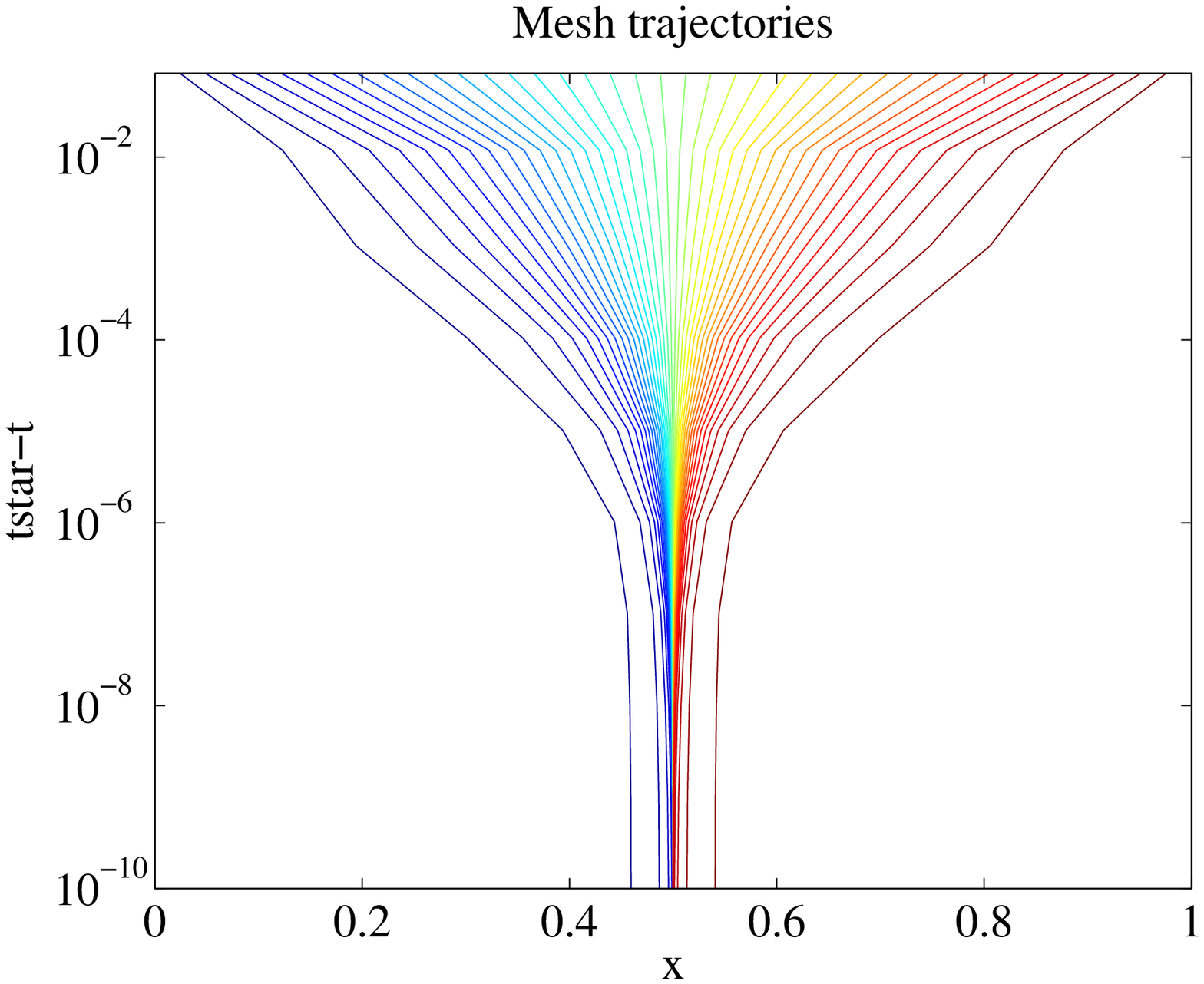}
  \caption{Plot of the solution profiles (left) and mesh contours (right)
    for the blow-up problem with $N=200$ and variable $\tau$ in the case
    $p=2$. The self-similar profile is displayed as a dashed line for
    comparison.}
  \label{fig:p2var}
\end{figure}

We next consider the effect of varying the mesh relaxation parameter
$\tau$ by selecting two constant values ($\tau=10^{-1}$ and $10^{-5}$)
as well as varying $\tau$ according to our strategy outlined in
Section~\ref{sec:varytau}.  The variable $\tau$ results are presented
for comparison purposes in Figure~\ref{fig:p2var}.  There is clearly
some loss of self-similarity in the solution profile relative to the
constant $\tau$ simulations, which leads to considerably more mesh
points being located outside the blow-up peak; however, the peak is
still reasonably well-resolved, and there are indeed a number of other
reasons that the variable-$\tau$ results are superior, which we discuss
next.

\begin{figure}[htbp]
  \centering
  \psfragfontsize
  \psfrag{tstar}[Bc][c]{$\tstar$}
  \psfrag{N}[c][b]{$N$}
  \psfrag{CPU secs.}[Bc][c]{\emph{CPU secs.}}
  \includegraphics[width=0.4\textwidth,clip]{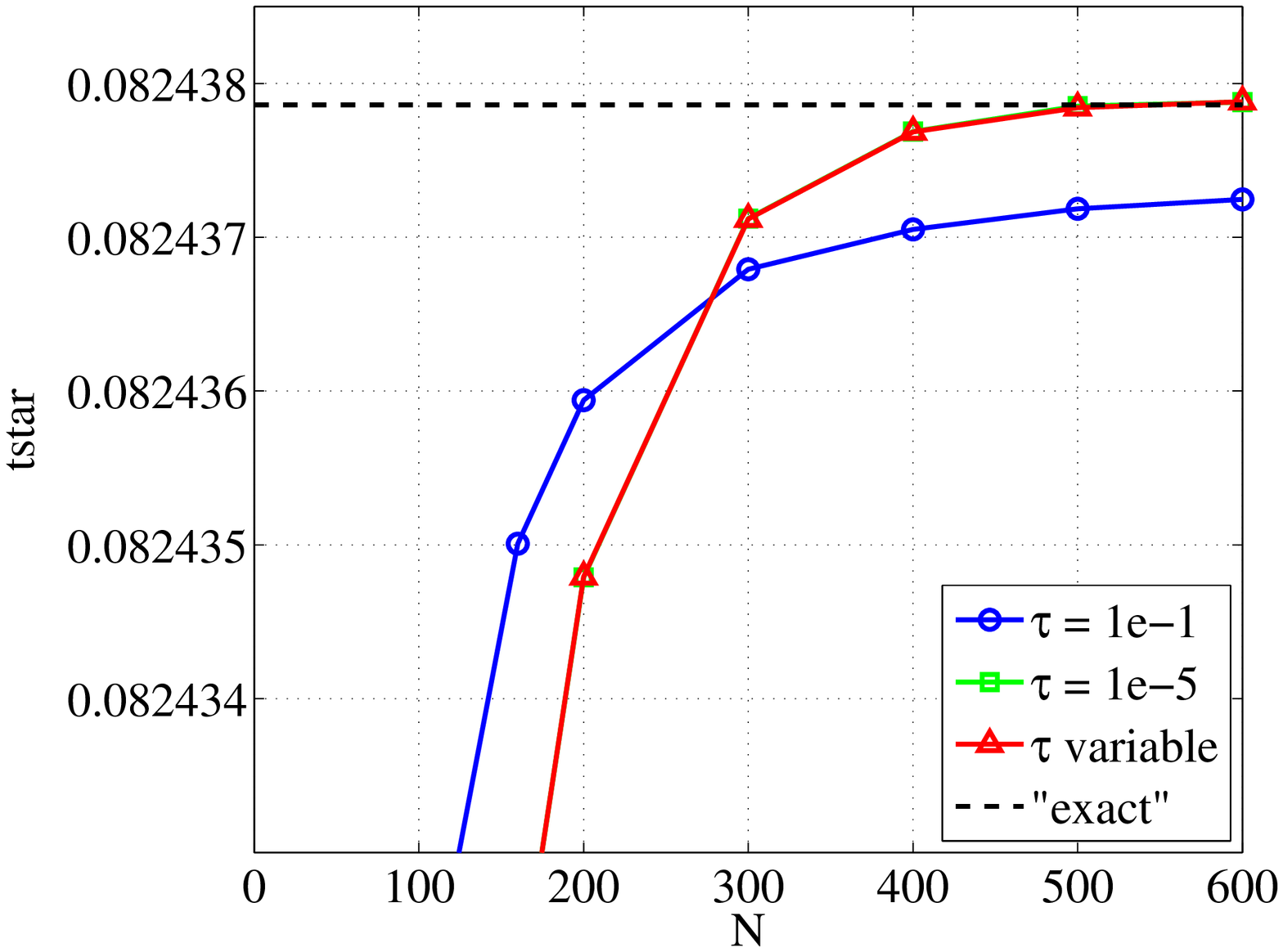}
  \quad
  \includegraphics[width=0.4\textwidth,clip]{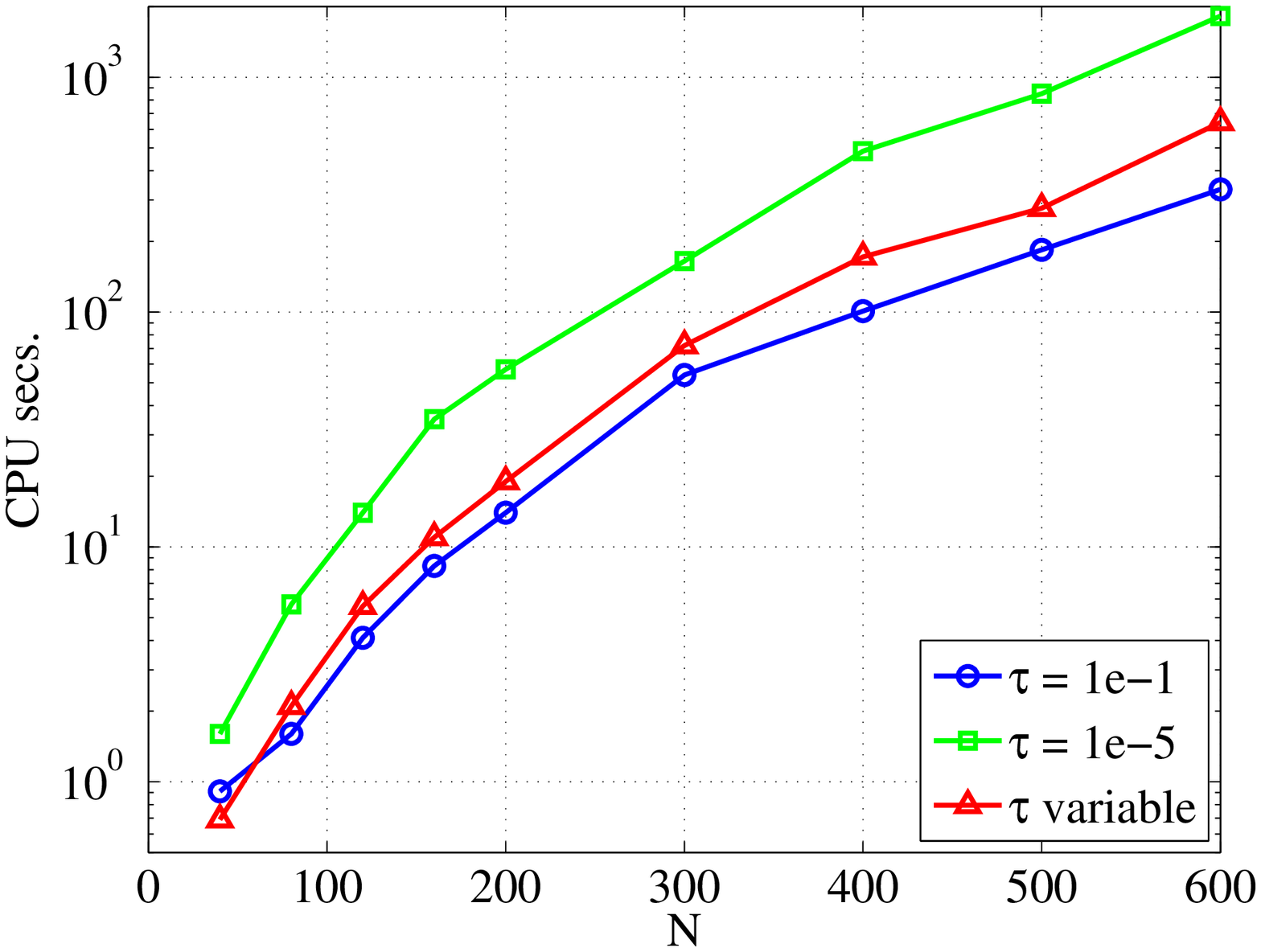}
  \caption{Comparison of blow-up time estimates (left) and CPU times
    (right) for various choices of $\tau$ in the case $p=2$.  The best
    estimate of the exact value of $\tstar \approx 0.08243786$ is shown
    as a dashed line.}
    %%The plot on the right is a zoomed-in view.}
  \label{fig:p2tstar}
\end{figure}

First, we performed a grid refinement study by
varying $N$ between 40 and 600, and compared the estimated blow-up times
for all choices of $\tau$ in Figure~\ref{fig:p2tstar}.  First of all,
the constant $\tau=10^{-5}$ result with $N=40$ points is consistent with
the value of $t^\star=0.082283$ reported in \cite{bhr96}.  The
$\tau=10^{-1}$ results require the least CPU time because such excessive
temporal smoothing acts to reduce the stiffness in the mesh equation;
however, the estimate of $\tstar$ is much less accurate and does not
converge to the correct blow-up time as the other simulations do.  Among
the remaining results ($\tau=10^{-5}$ and $\tau$ variable), there is no
visible difference in the blow-up time seems to suggest no advantage in
terms of accuracy.
Nonetheless, the variable $\tau$ approach is still capable of computing
further into the blow-up peak as evidenced by the maximum solution value
$\textsub{u}{max}$: for $\tau=10^{-5}$, all values of
$\textsub{u}{max}$ lie between $10^{14}$ and $10^{15}$ while for the
variable $\tau$ computations $\textsub{u}{max}$ is always above
$10^{18}$ at the end time.  There is only a slight loss of
self-similarity in the variable $\tau$ calculation, which can be seen by
comparing the solution with the asymptotic profile \en{ucos}, displayed
as a dashed line in Figures~\ref{fig:p2const}--\ref{fig:p2var}.

\begin{figure}[htbp]
  \centering
  \psfragfontsize
  \psfrag{tau}[Bc][c]{$\tau$}
  \psfrag{tstar-t}[c][b]{$\tstar-t$}
  \psfrag{Regularization parameter, tau}{}
  \includegraphics[width=0.4\textwidth,clip]{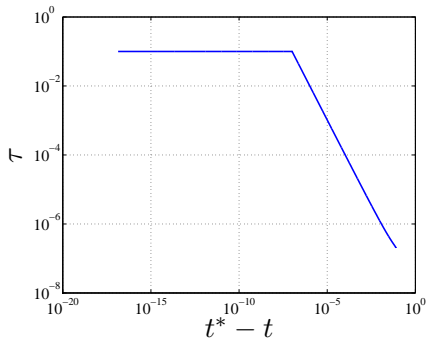}
  \caption{Plot of the mesh relaxation time, $\tau$, for the same
    parameters as for the $p=2$ problem with variable $\tau$.} 
  \label{fig:p2var-tau}
\end{figure}

The primary advantage of the variable $\tau$ approach is in terms of
efficiency, owing to the enhancement in temporal smoothing that occurs
close to the blow-up time.  The CPU time required in the variable $\tau$
case is consistently smaller by at least a factor of three relative to
the constant $\tau$ computations, as depicted in
Figure~\ref{fig:p2tstar}.  The variation of $\tau$ with time is shown in
Figure~\ref{fig:p2var-tau}, which demonstrates a nearly linear
dependence as the blow-up point is approached.  As a result, we can
think of the variable $\tau$ algorithm as keeping the mesh relaxation
time small when it is most needed (at the time when the blow-up peak is
first forming), but then introducing significant temporal smooting
closer to the blow-up time when the mesh equation is most stiff, even
though the mesh points themselves are not moving appreciably.

In summary, the use of a variable $\tau$ permits a more accurate
computation of both the blow-up time and the solution evolution, while
still maintaining a reasonable degree of self-similarity in the
solution, and all this at a significant savings in computational cost.
The primary reason for the improvement in performance is the reduction
in stiffness of the moving mesh PDE which results from allowing the mesh
relaxation time to vary.

\subsection{Blow-up with $p=5$}
\label{sec:results-p5}

\begin{figure}[htbp]
  \centering
  \psfragfontsize
  \psfrag{xi}[c][b]{$\xi$}
  \psfrag{\(u/umax\)}[Bc][c]{$(u/\textsub{u}{max})^4$}
  \psfrag{4}{}
  \psfrag{tstar-t}[c][b]{$\tstar-t$}
  \psfrag{x}[Bc][c]{$x$}
  \psfrag{Mesh trajectories}{}
  \includegraphics[width=0.4\textwidth,clip]{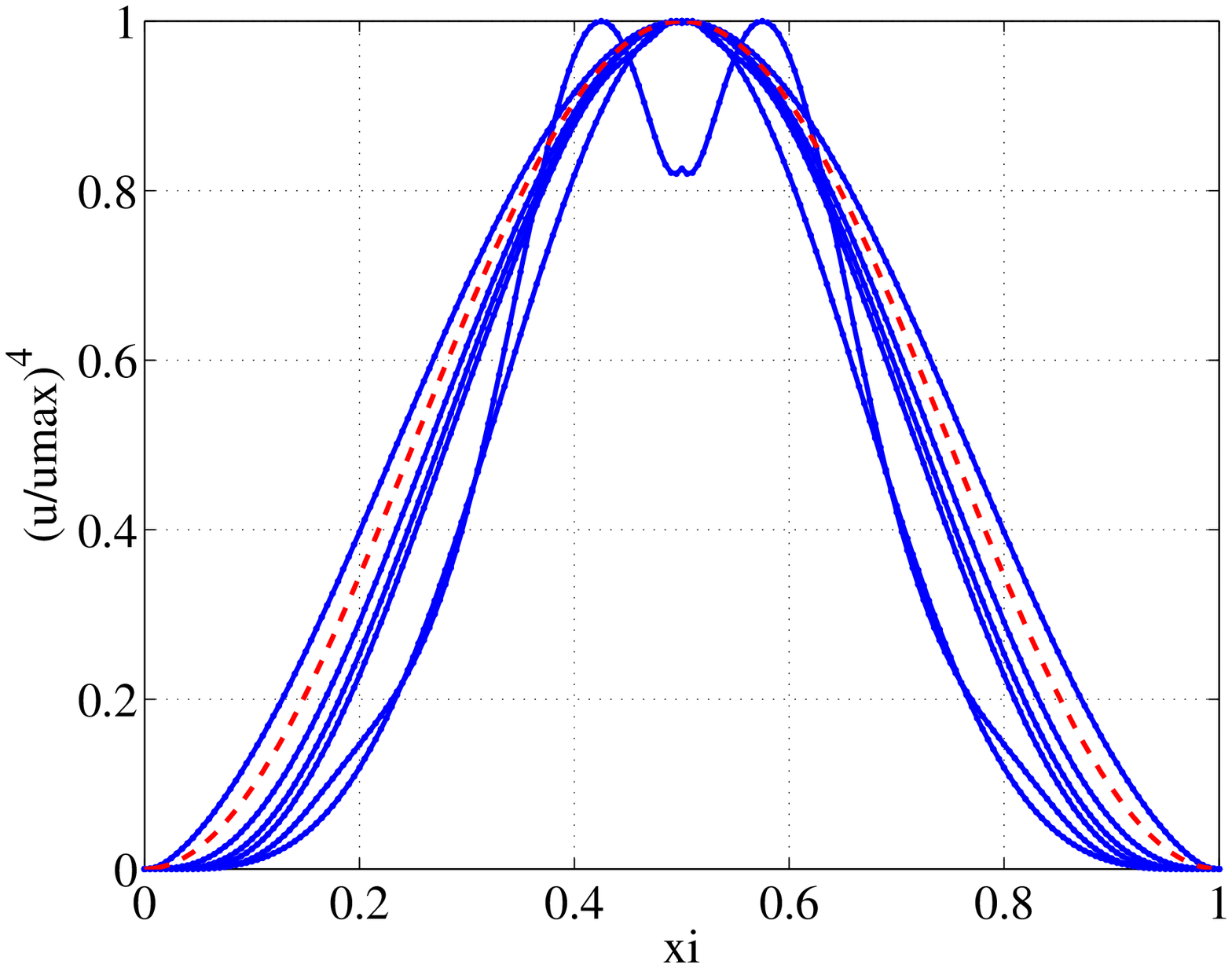} 
  \quad
  \includegraphics[width=0.4\textwidth,clip]{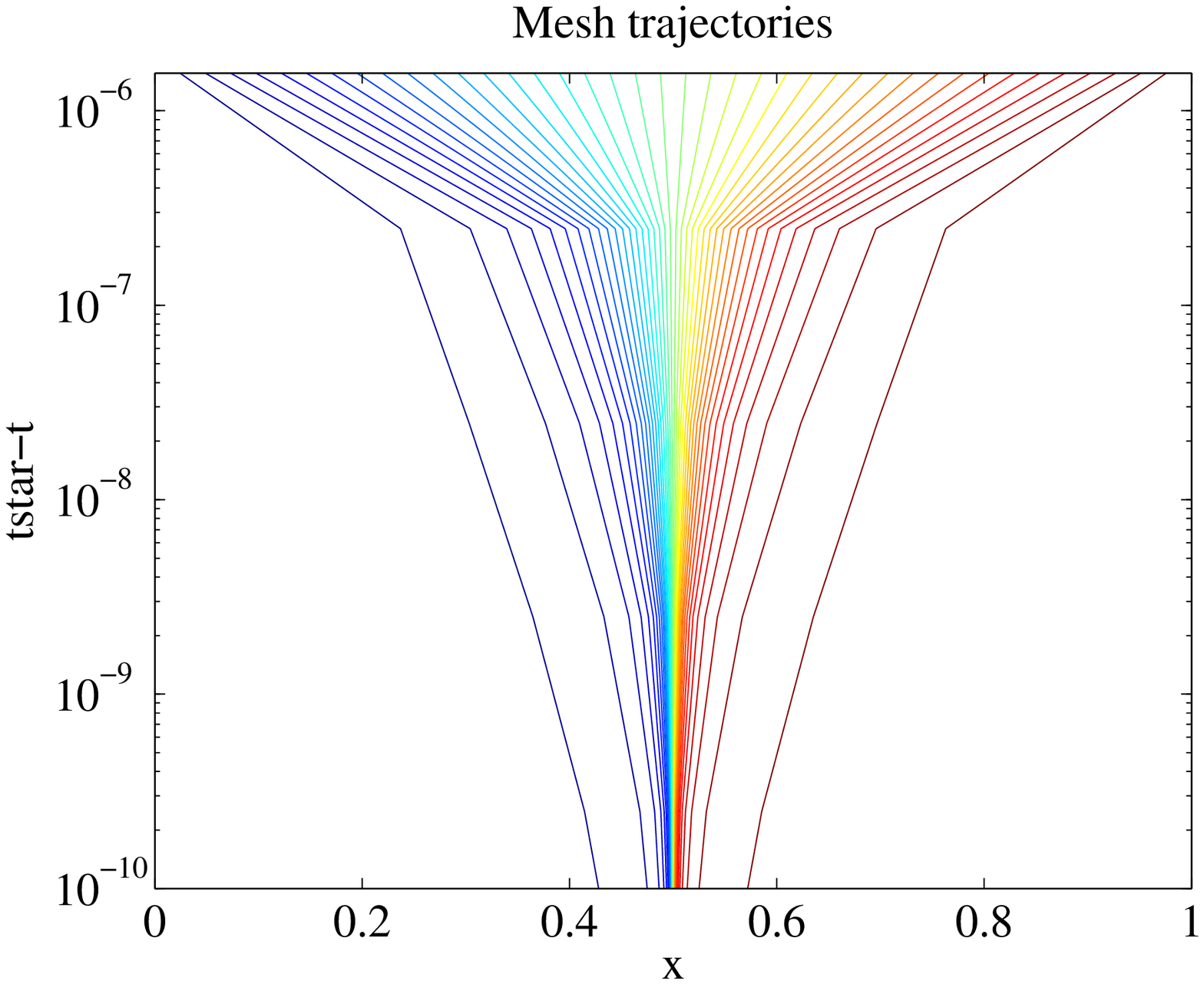} 
  \caption{Plot of the $p=5$ blow-up solution profiles (left) and mesh
    contours (right) for fixed $\tau=10^{-5}$ with $N=200$.  The
    self-similar profile is displayed as a dashed line for comparison.} 
  \label{fig:p5const}
\end{figure}

\begin{figure}[htbp]
  \centering
  \psfragfontsize
  \psfrag{xi}[c][b]{$\xi$}
  \psfrag{\(u/umax\)}[Bc][c]{$(u/\textsub{u}{max})^4$}
  \psfrag{4}{}
  \psfrag{tstar-t}[c][b]{$\tstar-t$}
  \psfrag{x}[Bc][c]{$x$}
  \psfrag{Mesh trajectories}{}
  \includegraphics[width=0.4\textwidth,clip]{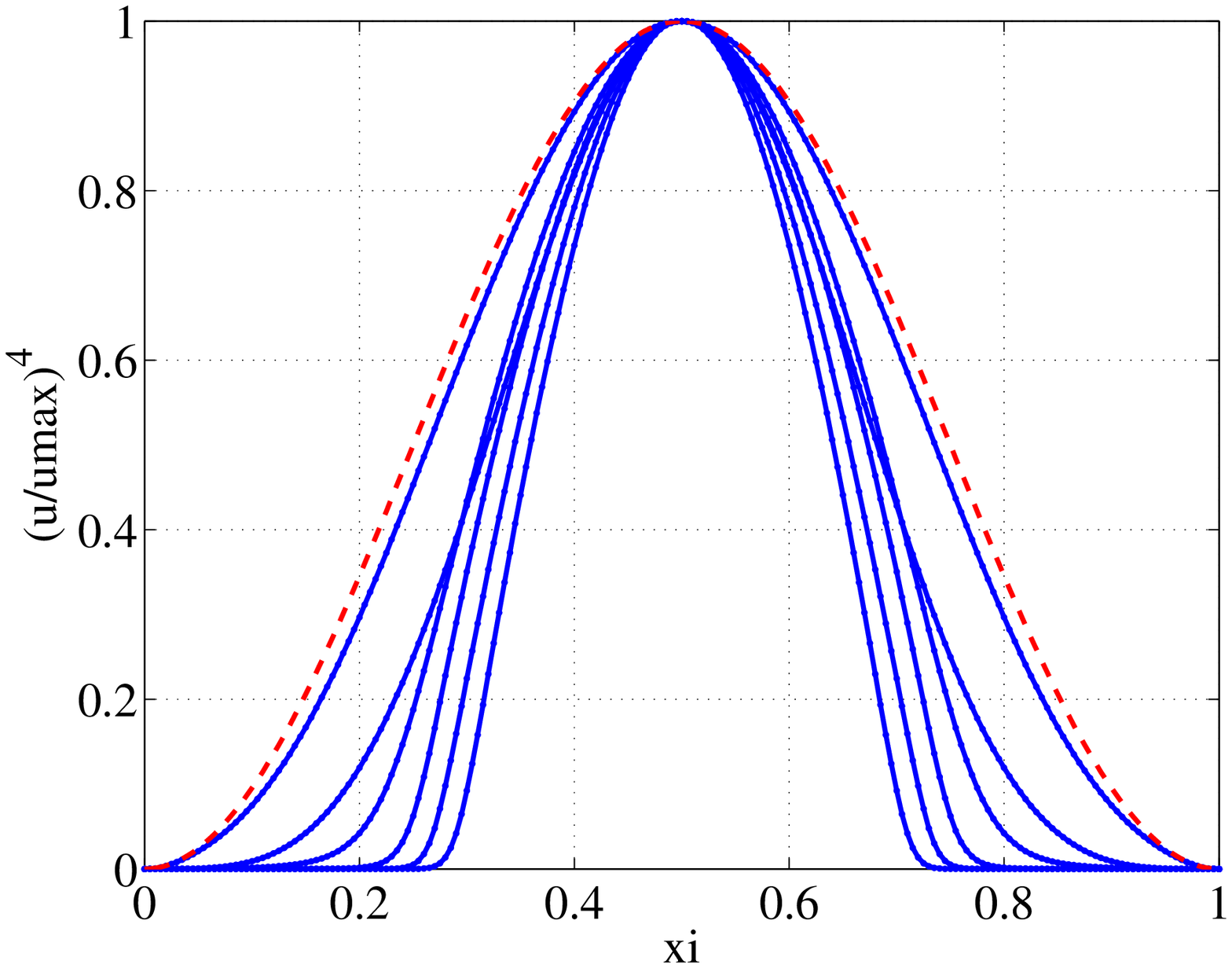}
  \quad
  \includegraphics[width=0.4\textwidth,clip]{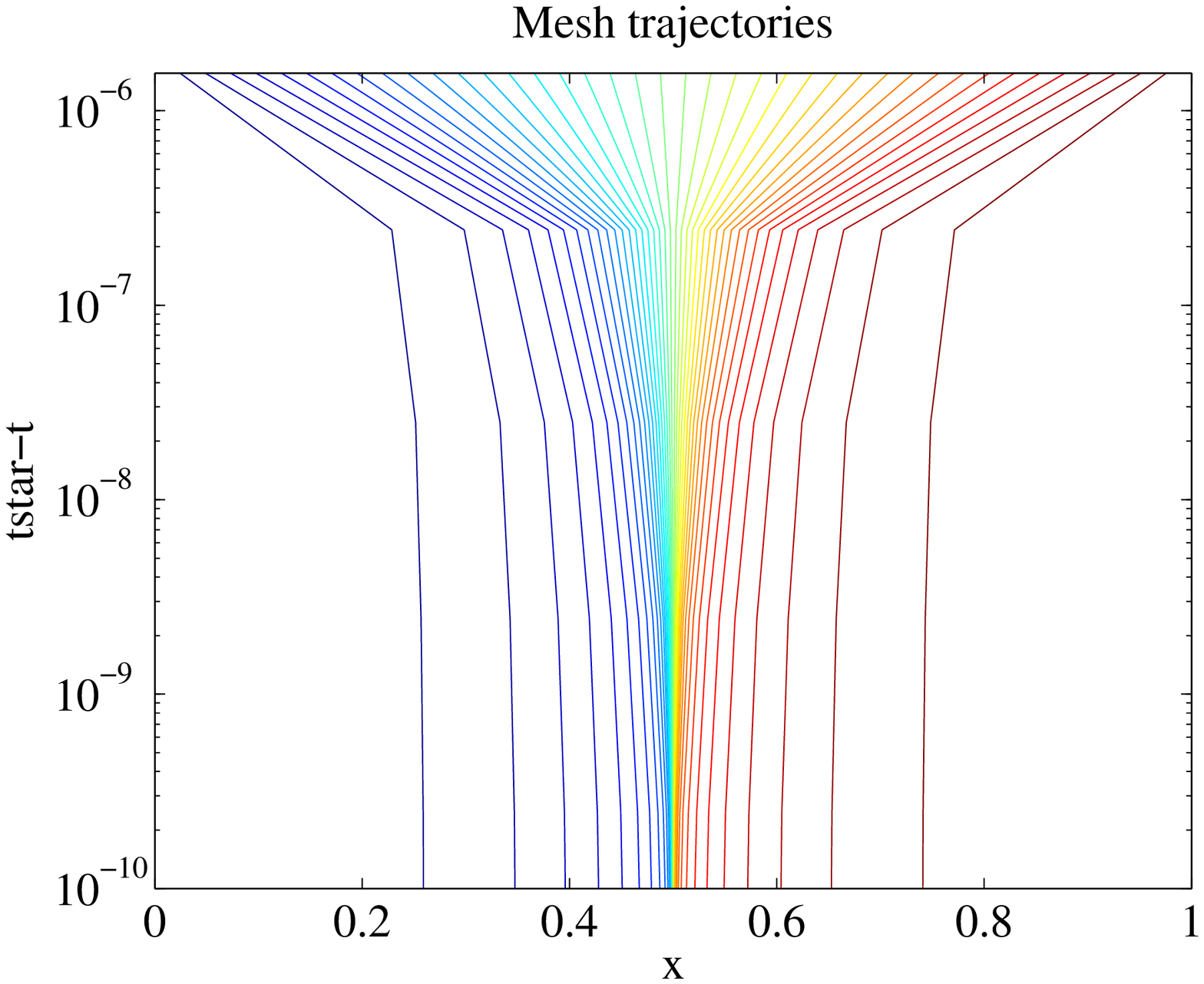}
  \caption{Plot of the $p=5$ blow-up solution profiles (left) and mesh
    contours (right) for variable $\tau$ in the case $N=200$.  The
    self-similar profile is displayed as a dashed line for comparison.}
  \label{fig:p5var}
\end{figure}

\begin{figure}[htbp]
  \centering
  \psfragfontsize
  \psfrag{tstar}[Bc][c]{$\tstar$}
  \psfrag{N}[c][b]{$N$}
  \psfrag{CPU secs.}[Bc][c]{\emph{CPU secs.}}
  \includegraphics[width=0.4\textwidth,clip]{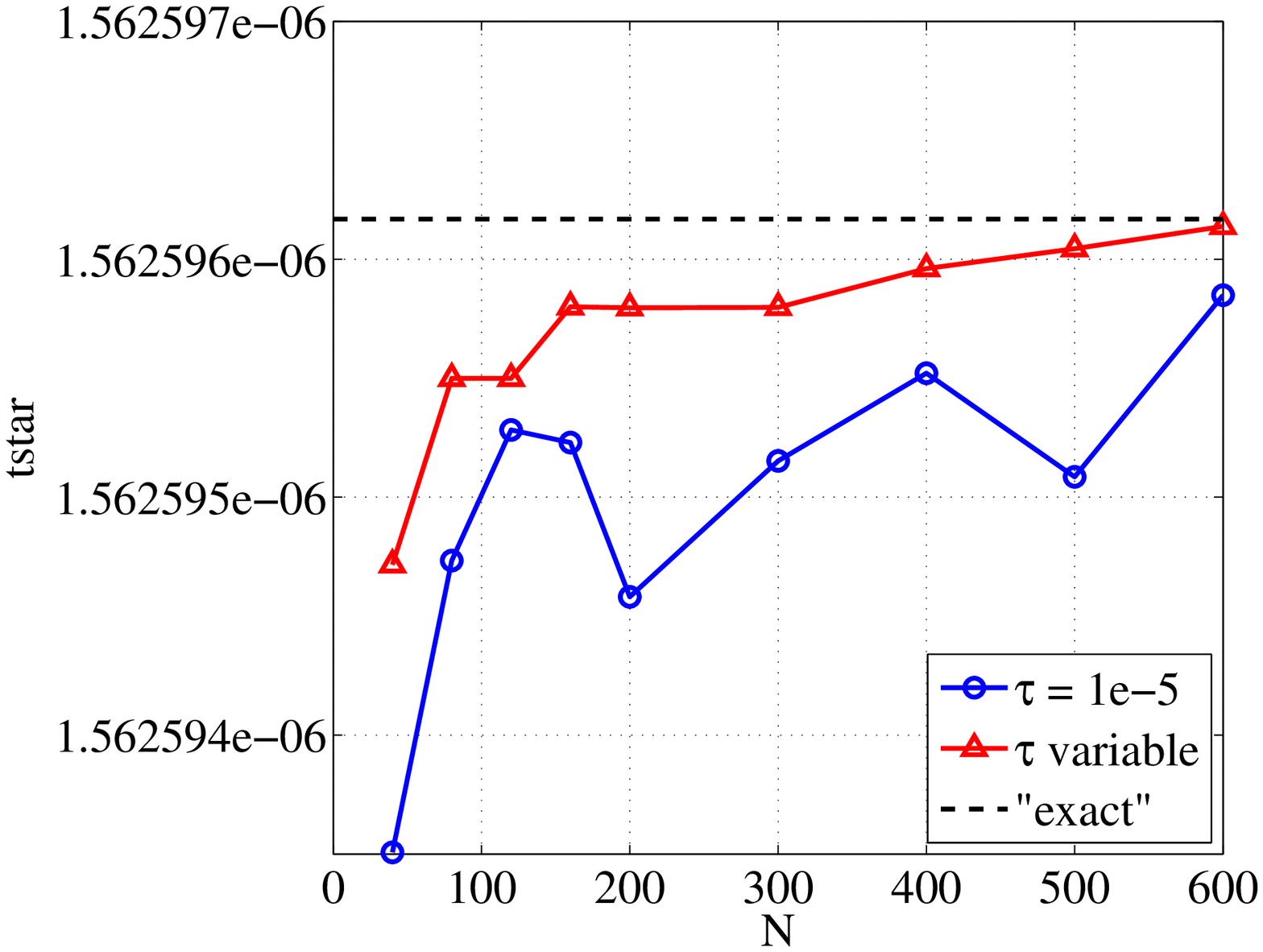}
  \quad 
  \includegraphics[width=0.4\textwidth,clip]{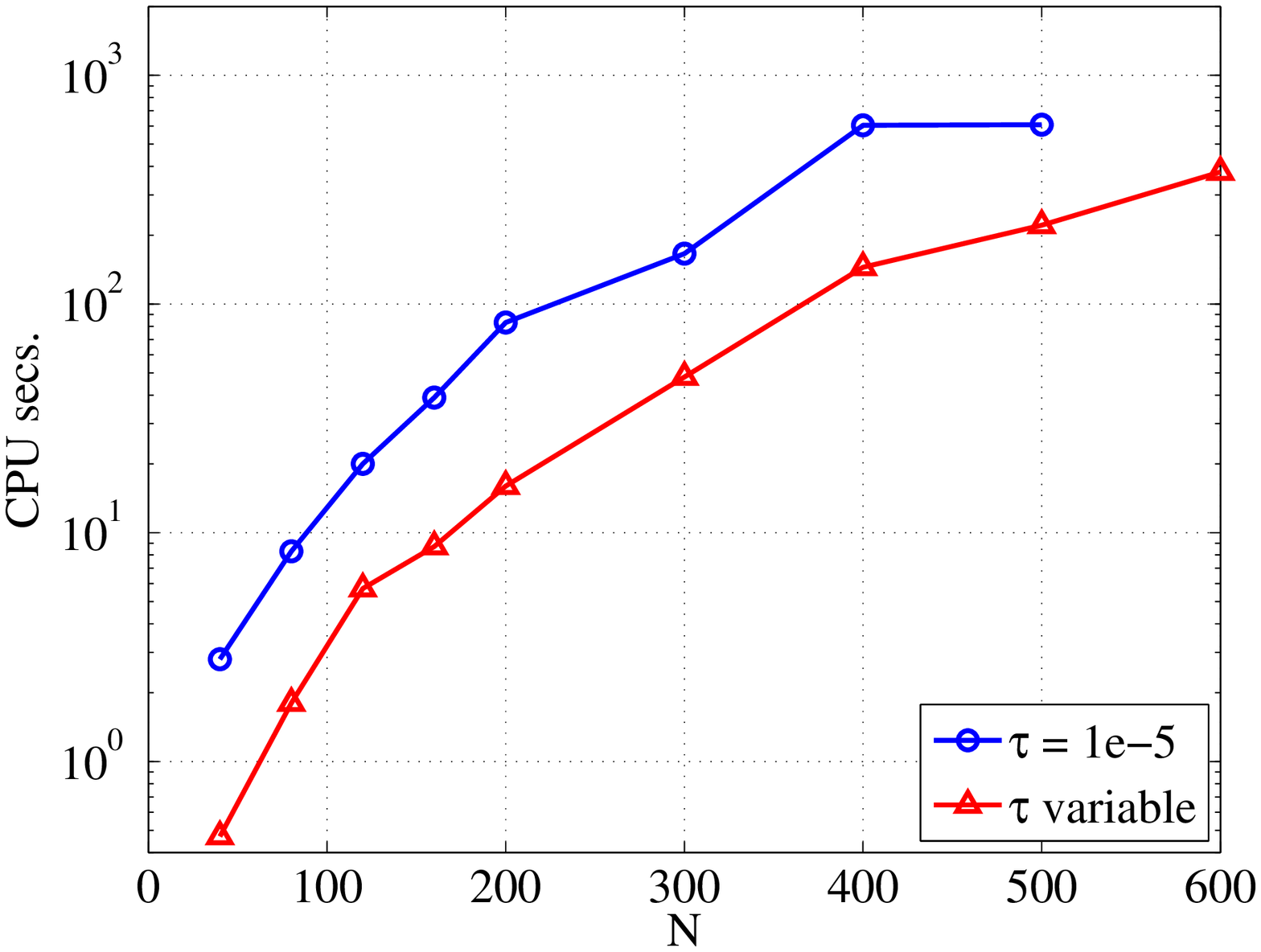}
  \caption{Comparison of blow-up time estimates
    (left) and CPU times (right) for the case $p=5$.  The best estimate
    of $\tstar \approx 1.5625962\times 10^{-6}$ (computed with $N=1000$
    points and increased DDASSL\ tolerances) is shown as a dashed
    line.}
  \label{fig:p5tstar}
\end{figure}

The $p=5$ blow-up problem constitutes a more difficult computational
problem, and so we consider it a more stringent test of our moving mesh
approach.  In this case, as in~\cite{bhr96}, we had to introduce mesh
smoothing ($\gamma=2$, $ip=4$) in order to ensure stability of the mesh
equation.  Proceeding as we did in the previous section, we compare the
$\tau=10^{-5}$ results to those for variable $\tau$, and the results are
depicted in Figures~\ref{fig:p5const} and \ref{fig:p5var}.
The constant $\tau$ computation exhibits oscillations in the solution
which cause the integration to fail due to numerical instability.  The
variable-$\tau$ results, on the other hand, show no such instability,
although the deviation from self-similarity is more significant than in
the $p=2$ case.  Nonetheless, the mesh points are still reasonably
well-clustered within the blow-up peak.  There is a similar three-fold
improvement in efficiency with the variable $\tau$ approach (see
Figure~\ref{fig:p5tstar}) although in this case the CPU times aren't as
meaningful because the constant $\tau$ computations fail owing to
mesh instability.  Again, we see the superiority of our adaptive
approach for solving blow-up problems.

The estimated blow-up times are plotted in Figure~\ref{fig:p5tstar},
which demonstrate further the instability experienced with the constant
$\tau$ computations.  The blow-up time for the variable $\tau$ result
converges nearly monotonically and we claim it is a much more accurate
estimate of the actual blow-up time for the $p=5$ calculation.

%%%%%%%%%%%%%%%%%%%%%%%%%%%%%%%%%%
\subsection{Exponential blow-up}
\label{sec:results-expu}

The following blow-up problem with an exponential nonlinearity
\begin{gather}
  u_t = u_{xx} + e^u,
  \label{eq:blowup-exp}
\end{gather}
was also considered in \cite{bhr96} and is an even more difficult test
of the moving mesh method.  The appropriate monitor function to use in
this case is $M(u)=e^u$, and in analogy with the derivation of
\en{ucos}, there exists an asymptotically self-similar profile
\begin{gather*}
  e^{(u-\textsub{u}{max})} \sim 
  \cos^2 \left( \pi(\xi - \textstyle{\frac{1}{2}}) \right). 
\end{gather*}
We start with initial data $u(x,0)=5\sin(\pi x)$, and use MMPDE6 with
spatial smoothing parameter $ip=4$.   

\begin{figure}[htbp]
  \centering
  \psfragfontsize
  \psfrag{xi}[c][b]{$\xi$}
  \psfrag{exp\(u-umax\)}[Bc][c]{$e^{(u-\textsub{u}{max})}$}
  \psfrag{4}{}
  \psfrag{tstar-t}[c][b]{$\tstar-t$}
  \psfrag{x}[Bc][c]{$x$}
  \psfrag{Mesh trajectories}{}
  \includegraphics[width=0.4\textwidth,clip]{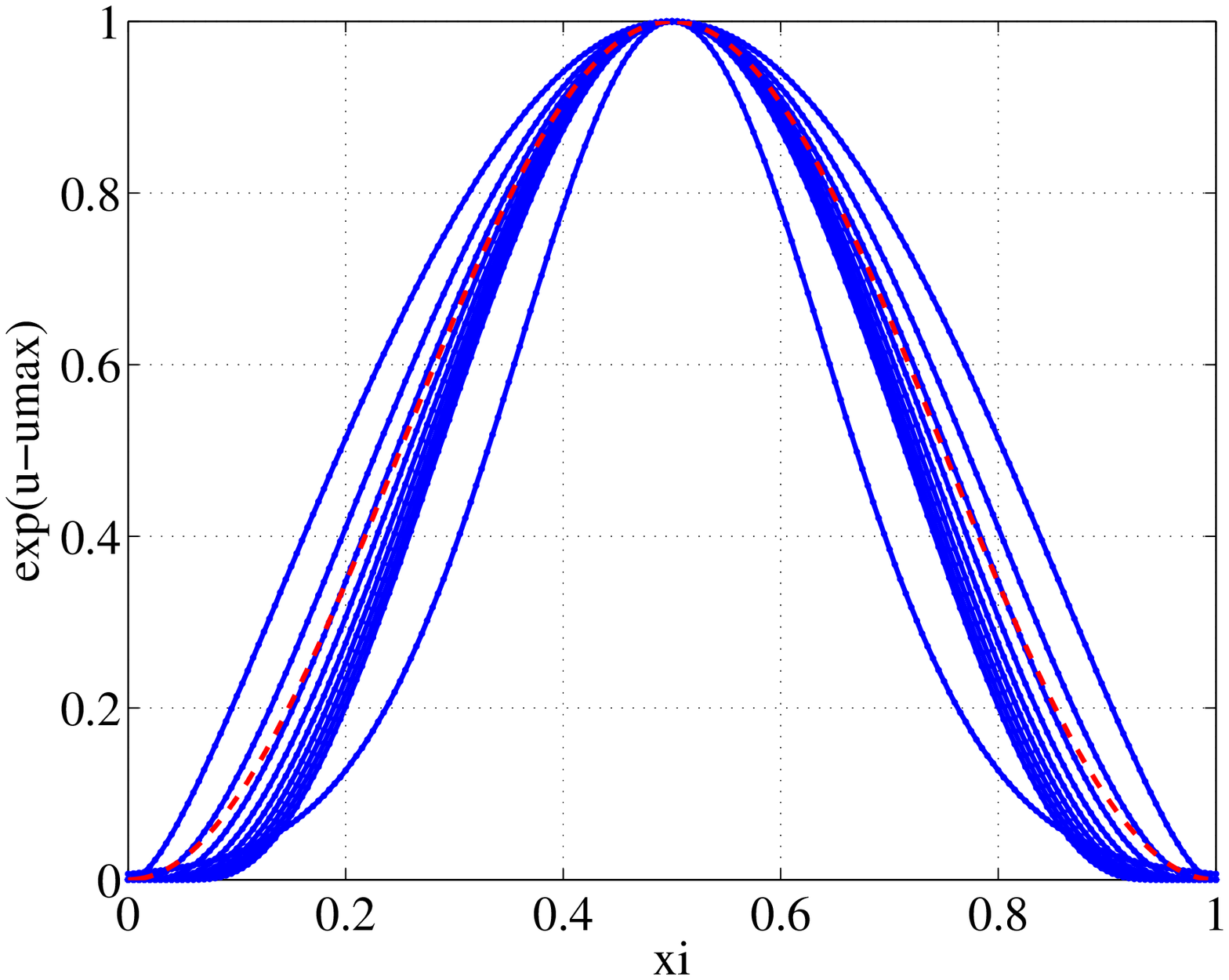} 
  \quad 
  \includegraphics[width=0.4\textwidth,clip]{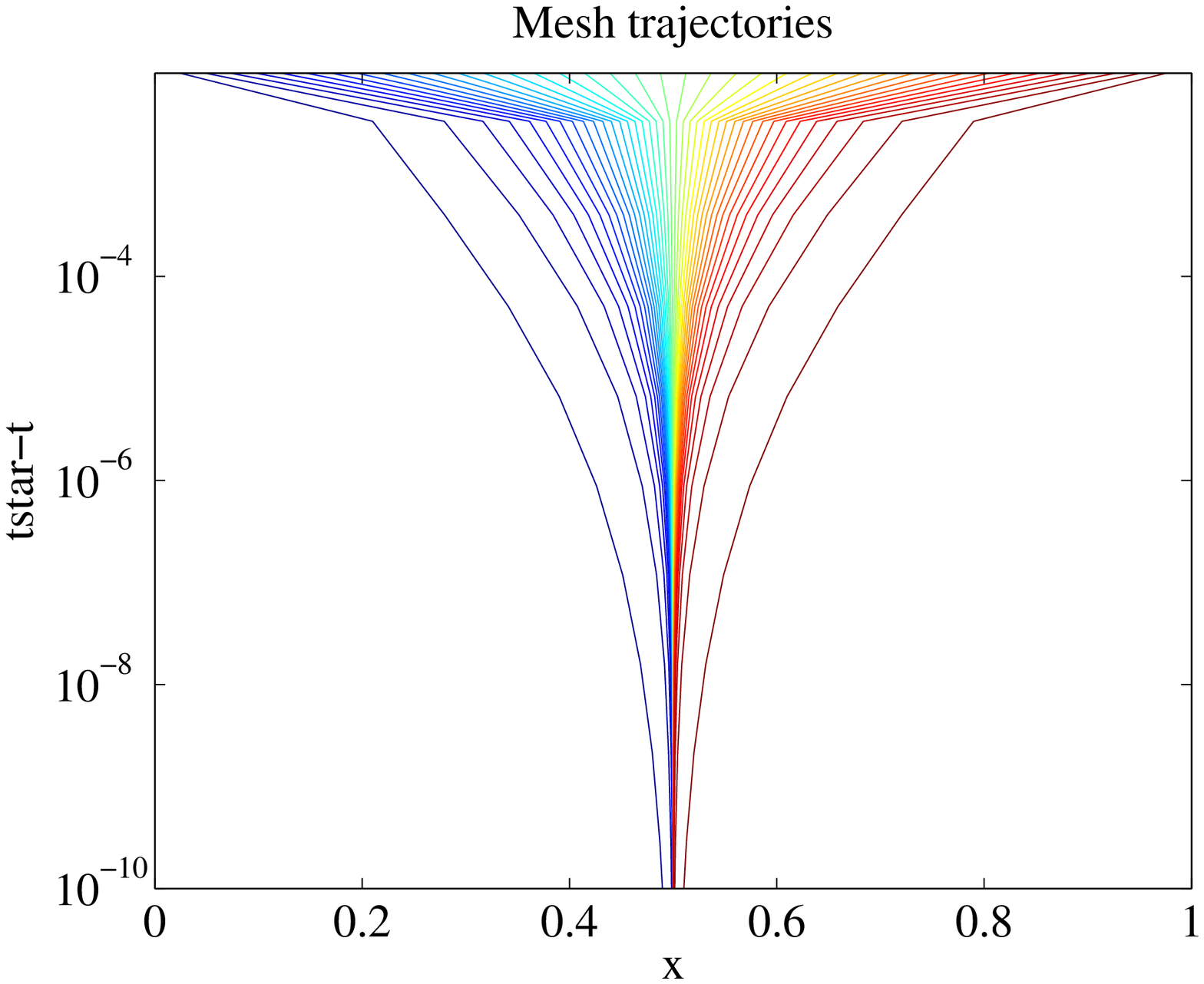} 
  \caption{Plot of the exponential blow-up solution profiles (left) and mesh
    contours (right) for fixed $\tau=10^{-5}$ with $N=200$.  The
    self-similar profile is displayed as a dashed line for comparison.} 
  \label{fig:euconst}
\end{figure}

\begin{figure}[htbp]
  \centering
  \psfragfontsize
  \psfrag{xi}[c][b]{$\xi$}
  \psfrag{exp\(u-umax\)}[Bc][c]{$e^{(u-\textsub{u}{max})}$}
  \psfrag{4}{}
  \psfrag{tstar-t}[c][b]{$\tstar-t$}
  \psfrag{x}[Bc][c]{$x$}
  \psfrag{Mesh trajectories}{}
  \includegraphics[width=0.4\textwidth,clip]{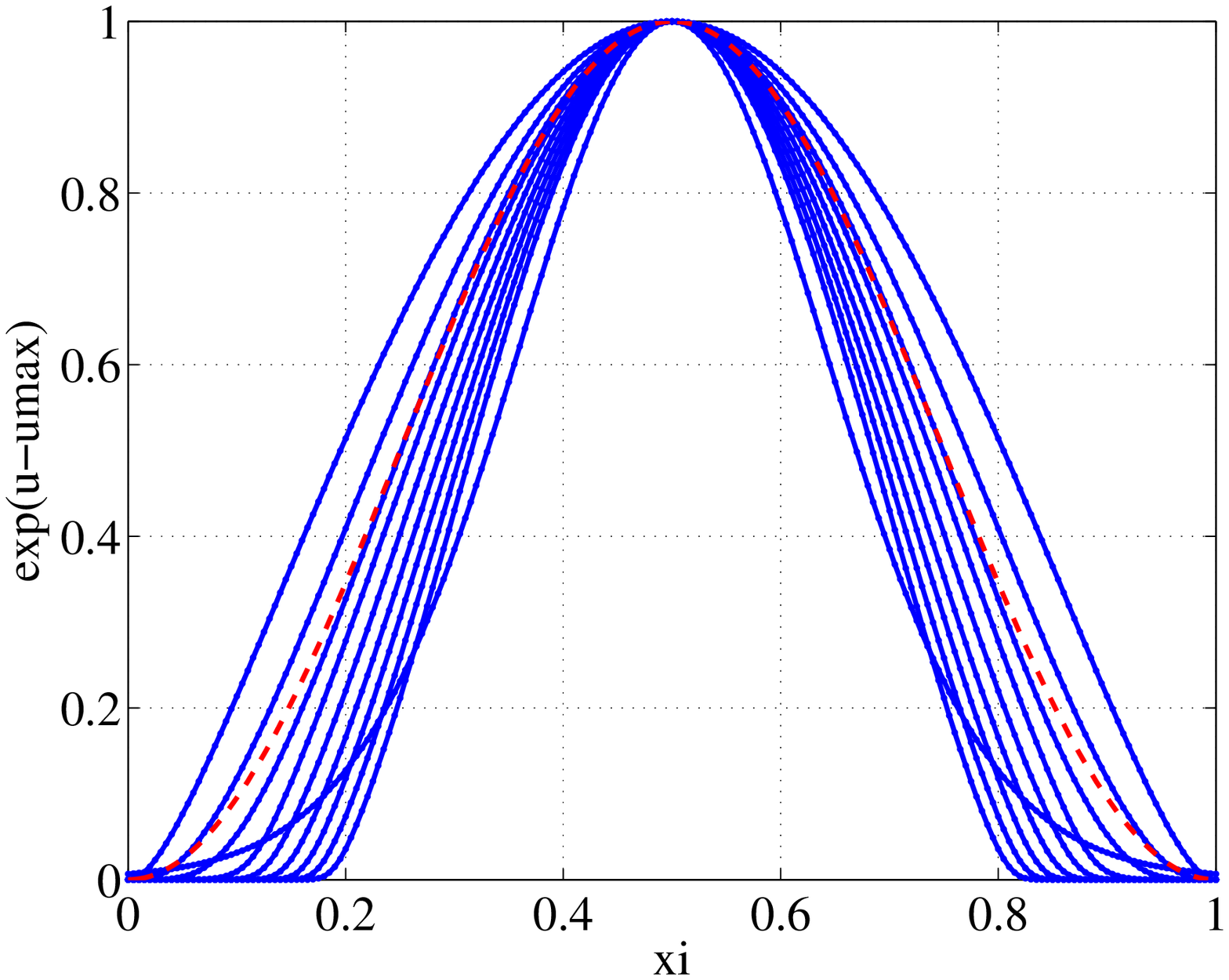}
  \quad 
  \includegraphics[width=0.4\textwidth,clip]{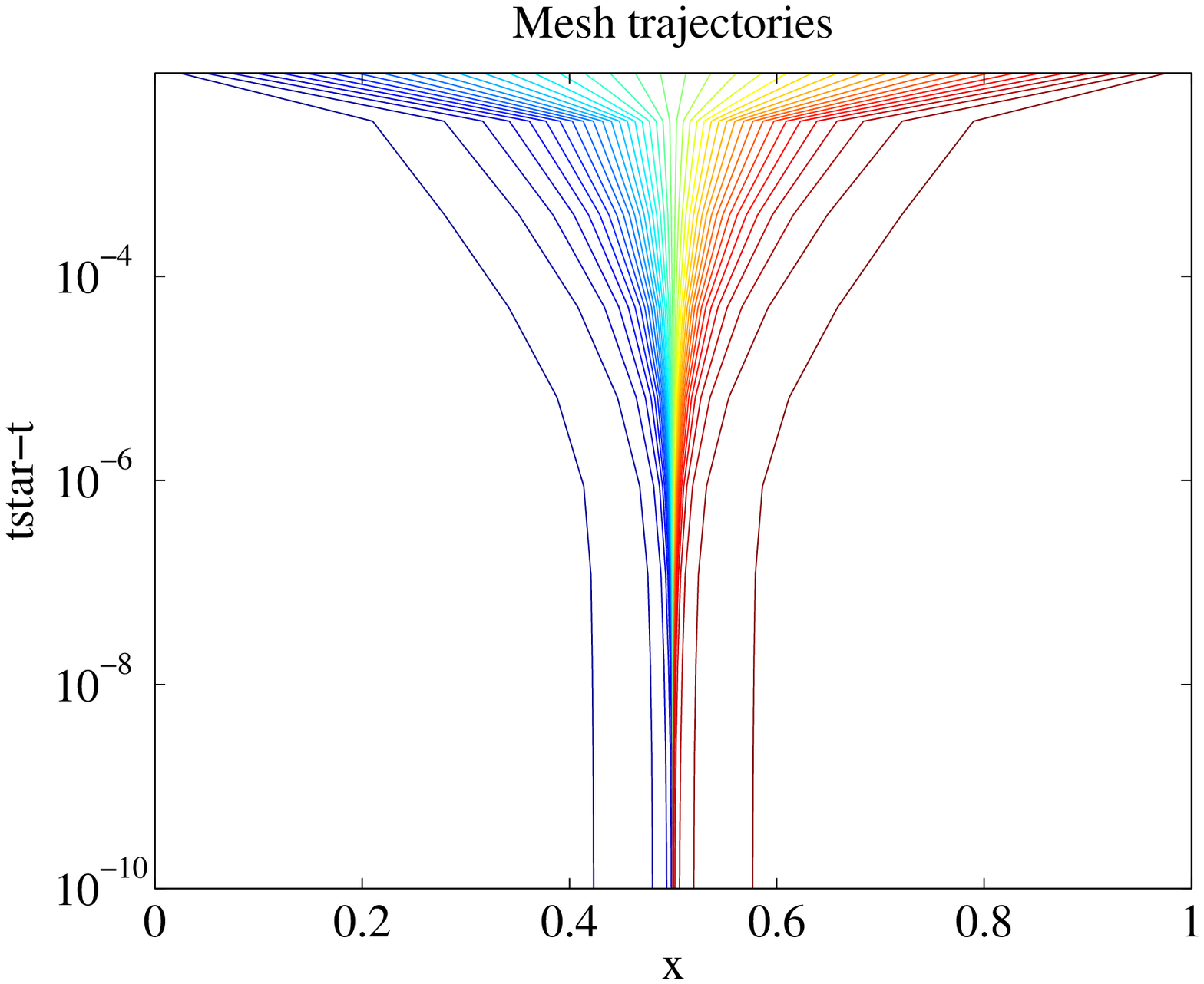}
  \caption{Plot of the exponential blow-up solution profiles (left) and mesh
    contours (right) for variable $\tau$ in the case $N=200$.  The
    self-similar profile is displayed as a dashed line for comparison.}
  \label{fig:euvar}
\end{figure}

\begin{figure}[htbp]
  \centering
  \psfragfontsize
  \psfrag{N}[c][b]{$N$}
  \psfrag{CPU secs.}[Bc][c]{\emph{CPU secs.}}
  \includegraphics[width=0.4\textwidth,clip]{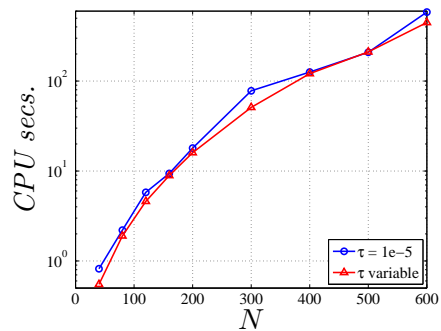}
  \caption{Comparison of CPU times (right) for the exponential blow-up
    problem.}
  \label{fig:eucpu}
\end{figure}

The results for $N=200$ mesh points are displayed in
Figures~\ref{fig:euconst} and~\ref{fig:euvar}.  There is a slight loss
of self-similarity in both cases owing to the introduction of spatial
smoothing, but the difference between the two solutions is minimal.  The
primary difference is in terms of efficiency, where the variable-$\tau$
simulation requires consistently 30\%\ less CPU time than for fixed
$\tau$.  This is not as dramatic an improvement as for the polynomial
blow-up examples considered in the previous two sections, as can be seen
in Figure~\ref{fig:eucpu}, but it is still a significant improvement.

%%%%%%%%%%%%%%%%%%%%%%%%%%%%%%%%%%
\leavethisout{
  \subsection{Gierer-Meinhardt and Gray-Scott Models}
  \label{sec:results-gmgs}
  
  \begin{note}
    Wentao Sun insists that this problem is too hard for moving meshes.
  \end{note}
  
  The Gierer-Meinhardt and Gray-Scott models consist of two coupled
  nonlinear reaction-diffusion equations which don't exhibit blow-up but
  rather propagating \emph{spikes} with the following types of behaviour:
  \begin{itemize}
  \item nearly stationary spikes which suddenly begin moving,
  \item two spikes which coalesce,
  \item a single spike for which the height oscillates,
  \item variable spike speed where the position oscillates in time. 
  \end{itemize}
  All of these behaviours arise due to bifurcations in the eigenvalues of
  the system and are dealt with in some detail in David Iron's Ph.D.
  thesis\cite{iron01}.  Two problems suggested by David
  are based on the following coupled system of two PDEs:
  \begin{gather}
    a_t = \epsilon_1^2 a_{xx} - a + \frac{a^2}{h}, \label{eq:gierer-meinhardt1}\\
    \epsilon_2 h_t = D h_{xx} - h + \frac{a^2}{\epsilon_1}, \label{eq:gierer-meinhardt2}
  \end{gather}
  where 
  \begin{gather*}
    a_x(\pm 1) = 0,
    h_x(\pm 1) = 0.
  \end{gather*}
  The two parameters $\epsilon_1$ and $\epsilon_2$ are typically small,
  and we will take $\epsilon_1=0.02$ and $\epsilon_2= 0.00001$.  We
  consider two cases:
  \begin{enumerate}
  \item $D=0.4$ and $a(x,0)=\frac{3}{2} \sech^2 \left(
      \frac{x+0.4}{2\epsilon_1} \right) + \frac{3}{2} 
    \sech^2 \left( \frac{x-0.401}{2\epsilon_1} \right)$: whose solution
    consists of two spikes that initially move very slowly, but then one
    spike spontaneously (and very rapidly) vanishes.
  \item $D=0.38$ and $a(x,0)=\frac{3}{2} \sech^2 \left(
      \frac{x+0.5}{2\epsilon_1} \right) + \frac{3}{2} 
    \sech^2 \left( \frac{x-0.5}{2\epsilon_1} \right)$: whose solution
    reaches a seeming equilibrium, but then very quickly goes unstable.
  \end{enumerate}
  Both of these problems are excellent test problems for our
  variable-$\tau$ approach, since they clearly involve two different time
  scales on which solution variations occur, and consequently are problems
  for which we would like our method to also adapt the mesh appropriately.
}

%%%%%%%%%%%%%%%%%%%%%%%%%%%%%%%%%%%%%%%%%%%%%%%%%%%%%%%%%%%%%%%%%%%%%%%%
\leavethisout{
  \subsection{Generalized KdV equation}
  
  JF Williams suggested the following problem, known as a \emph{generalized
    KdV equation}: 
  \begin{gather*}
    u_t  = \pm(u_{xx} + u^p)_x,
  \end{gather*}
  which has a travelling wave solution
  \begin{gather*}
    u=u(x-ct) = \left[ \frac{c(p+1)}{2} \sech \left(
        \frac{\sqrt{c}(p-1)}{2}(x-ct) \right) \right]^{1/(p-1)}
  \end{gather*}
  This problem is stable only for $p<5$, and blow-up occurs for $p>5$.  
}

%%%%%%%%%%%%%%%%%%%%%%%%%%%%%%%%%%%%%%%%%%%%%%%%%%%%%%%%%%%%%%%%%%%%%%%%
\section{Conclusion}
\label{sec:conclusion}

In this paper, we have considered a moving mesh approach for solving
self-similar blow-up problems.  The novelty of this method stems from
its use of a solution-dependent mesh relaxation time, $\tau$.  We have
proposed a strategy for selecting $\tau$ in the context of self-similar
blow-up problems.  Numerical simulations demonstrate that by varying the
relaxation time in an appropriate way, the solution can be computed more
accurately, further into the blow-up, and more efficiently than
would otherwise be possible with a constant value of $\tau$.

Because our strategy for adapting $\tau$ is specific to problems of
blow-up type, we plan in the future to extend these results by
generalizing them to a more generic class of problems.  We intend to
investigate other nonlinear parabolic problems that exhibit more general
blow-up behaviour (such as the generalized Korteweg-de Vries or
Gierer-Meinhardt equations) as well as problems with moving fronts.

%%%%%%%%%%%%%%%%%%%%%%%%%%%%%%%%%%%%%%%%%%%%%%%%%%%%%%%%%%%%%%%%%%%%%%%%
\leavethisout{
  \section*{Outstanding issues}
  
  Others (for example, \cite{blp01}) make use of a time scaling,
  known as the \emph{Sundman transformation} of the form $dt/ds=1/M$,
  to eliminate difficulties with time integration.  This may help us
  here but requires transforming equations to the $s$ variable and
  adding one further equation for $t(s)$.
}

%%%%%%%%%%%%%%%%%%%%%%%%%%%%%%%%%%%%%%%%%%%%%%%%%%%%%%%%%%%%%%%%%%%%%%%%
\appendix
\section{Analysis of moving mesh equation}
\label{sec:taylor}

We briefly redo the analysis from \cite{hrr94b} for
time-dependent $\tau(t)$.  We begin with \en{mmesh-tau} and
perform a Taylor series expansion for small $\tau$ to obtain:
\begin{align*}
  \frac{\partial}{\partial \xi} x(\xi, t+\tau(t)) =\,& 
  \frac{\partial}{\partial \xi} x(\xi,t) + \tau (1+\dot{\tau}) 
  \frac{\partial}{\partial \xi} \dot{x}(\xi,t) + 
  \order{\tau^2},\\
  \intertext{and}
  M(x(\xi, t + \tau(t)),t+\tau(t)) =\,& 
  M(x(\xi,t), t) + \\
  & \tau (1+\dot{\tau}) \left(
    \dot{x}\, \frac{\partial}{\partial \xi} M(x(\xi,t),t) 
    + \frac{\partial}{\partial t} M(x(\xi,t),t) \right)
  + \order{\tau^2}.
\end{align*}
Substituting these expressions into \en{mmesh-tau} we obtain the
corresponding equations for MMPDE4 and MMPDE6 respectively:
\begin{align} 
  \tau(1+\dot{\tau})\, \frac{\partial }{\partial \xi} 
  \left(M\frac{\partial \dot{x} }{\partial \xi}\right) &=
  - \frac{\partial }{\partial
    \xi}\left(M \frac{\partial x
  }{\partial \xi}\right), \\
  \tau(1+\dot{\tau})\, \frac{\partial ^2 \dot{x} }{\partial \xi^2} &= -
  \frac{\partial }{\partial \xi}(M \frac{\partial x}{\partial \xi}).
\end{align}
Notice that relative to \en{mmpde4} and \en{mmpde6}, the only
change here is an extra factor of $(1+\dot{\tau})$ which simply scales
$\tau$.  Therefore a time-dependent $\tau$ has a minimal impact on
the moving mesh equation.

%%%%%%%%%%%%%%%%%%%%%%%%%%%%%%%%%%%%%%%%%%%%%%%%%%%%%%%%%%%%%%%%%%%%%%%%
%%\begin{ack}
%%  This work was supported by a grant from the Natural Sciences and
%%  Engineering Research Council of Canada (NSERC).
%%\end{ack}

%%%%%%%%%%%%%%%%%%%%%%%%%%%%%%%%%%%%%%%%%%%%%%%%%%%%%%%%%%%%%%%%%%%%%%%%
\newcommand{\bibvol}[2]{#1}
\newcommand{\bibtitle}[1]{#1}
\newcommand{\bibjournal}[1]{\emph{#1\/}}

\bibliographystyle{plain}

\end{document}